\documentclass{article}
\usepackage{graphicx} % Required for inserting images
\usepackage[utf8]{inputenc}
\usepackage[margin=1in]{geometry}
\usepackage{amsmath}
\usepackage{amsthm}
\usepackage{amsfonts}
\usepackage{amssymb}
\usepackage{txfonts}
\usepackage[backend=biber, sorting=none]{biblatex}
\usepackage{ragged2e}
\usepackage{authblk}
\usepackage{hyperref}

\providecommand{\keywords}[1]
{
  \small	
  \textbf{\textit{Keywords: }} #1
}

\providecommand{\subjclass}[1]
{
  \small	
  \textbf{\textit{Mathematics Subject Classification: }} #1
}

\newtheorem{theorem}{Proposition}
\newtheorem{corollary}[theorem]{Corollary}
\title{Generalization of some of Ramanujan's formulae}
\author{\emph{Aung Phone Maw}}
\date{August 2024}
\begin{document}

\maketitle
\begin{abstract}
    We shall make use of the method of partial fractions to generalize some of Ramanujan's infinite series identities, including Ramanujan's famous formula for $\zeta(2n+1)$, and we shall also give a generalization of the transformation formula for the Dedekind eta function. It is shown here that the method of partial fractions can be used to obtain many similar identities of this kind.
\end{abstract}
\raggedright
. \\
.
\par
\section{Introduction}

The well-known Eulerian formula $ \zeta(2n) = \frac{(-1)^{n-1}(2\pi)^{2n}B_{2n}}{2(2n)!}$, where $B_{2n}$ are the Bernoulli numbers, provides an exact evaluation of the zeta function at even integers. On the other hand no closed-form expression for the values of the zeta function at odd integers are known. In that direction, Lerch [10] and Ramanujan [1, 8] provided impressive identities giving rapidly converging series for $\zeta(2n+1)$. Lerch formula is known to be just a particular case of Ramanujan's more general formula : 

\begin{multline*}
    -2^{2n+1} \sum\limits_{j=0}^{n+1} (-1)^j \alpha^{n+1-j} \beta^j \frac{B_{2j}B_{2n+2-2j}}{(2j)! (2n+2-2j)!} = \alpha^{-n}\left\{ \zeta(2n+1) + 2\sum\limits_{r=1}^{\infty} \frac{1}{r^{2n+1}(e^{2r \alpha}-1)} \right\}-  (-\beta)^{-n}\left\{ \zeta(2n+1) + 2\sum\limits_{r=1}^{\infty} \frac{1}{r^{2n+1}(e^{2r \beta}-1)} \right\}.
\end{multline*}
Which holds for $\alpha\beta = \pi^2$. Various proofs of this formula have been investigated and studied, via the usage of Mittag-Leffler expansion [11], residue calculus [14, 15], or the partial fraction expansion of coth function [12, 13]. A comprehensive overview is given by Berndt and Straub [3].
\\ \vspace{0.7cm}
In our investigation, we shall utilize and manipulate the elementary partial fraction decomposition relations to arrive at even more generalized expressions. We will arrive at generalizations of not only Ramanujan's expressions, we will also arrive at a generalization of the transformation formula for the Dedekind eta function $\eta(-\frac{1}{\tau})= (-i\tau)^{-\frac{1}{2}} \eta(\tau)$. Where $\eta(\tau) = e^{\frac{\pi i \tau}{12}}\prod\limits_{n=1}^{\infty} (1-e^{2\pi i n \tau})$, $\textbf{\textit{Im}}(\tau) > 0$. We shall now state some elementary partial fraction relations. It is well-known that the following partial fraction decomposition identities hold. Provided that each $x_i, 1 \leq i \leq M$ is distinct.

\[ \tag{1.1} \label{eqb1} \frac{1}{\prod\limits_{1\leq i \leq M}(x_i-t)}  =\sum\limits_{1 \leq i \leq M} \frac{1}{(x_i-t)\prod\limits_{\substack{1 \leq j \leq M \\ j \neq i }}(x_j-x_i)}.\]

\[ \tag{1.2} \label{eqab1} \frac{1}{\prod\limits_{1\leq i \leq M}(1-x_i t)}  =\sum\limits_{1 \leq i \leq M} \frac{1}{(1-x_i t)\prod\limits_{\substack{1 \leq j \leq M \\ j \neq i }}(1-\frac{x_j}{x_i})}.\]

Now for our purpose, let us make further adjustments to \eqref{eqab1}. Let us define another set of $N$ variables $y_1,...,y_N$ and in \eqref{eqab1} let us set $t = y_i$, where $ 1 \leq i \leq N$ and multiply both sides by $\frac{1}{\prod\limits_{\substack{1 \leq j \leq N \\ j \neq i }}(1-\frac{y_j}{y_i})}$ and then let us sum through $1 \leq i \leq N$. Then we have :

\[\tag{1.3} \label{eqab2}
\sum\limits_{1 \leq i \leq N}\frac{1}{\prod\limits_{1 \leq j \leq M}(1-x_jy_i)\prod\limits_{\substack{1 \leq j \leq N \\ j \neq i }}(1-\frac{y_j}{y_i})}= \sum\limits_{1 \leq i \leq M}\frac{1}{\prod\limits_{1 \leq j \leq N}(1-x_i y_j)\prod\limits_{\substack{1 \leq j \leq M \\ j \neq i }}(1-\frac{x_j}{x_i})}.
\]

Now let $N \to N+1$ and let us set $y_{N+1}$ to be $t$. Then we arrive to the following partial fraction identity.

\begin{multline}
    \tag{1.4} \label{eqab3}
 \frac{1}{\prod\limits_{1 \leq i \leq M}(1-x_i t)\prod\limits_{1 \leq i \leq N}(1-\frac{y_i}{t})} = \sum\limits_{1 \leq i \leq M}\frac{1}{(1-x_i t)\prod\limits_{1 \leq j \leq N}(1-x_iy_j)\prod\limits_{\substack{1 \leq j \leq M \\ j \neq i }}(1-\frac{x_j}{x_i})} \\ -\sum\limits_{1 \leq i \leq N}\frac{1}{(1-\frac{t}{y_i})\prod\limits_{1 \leq j \leq M}(1-x_j y_i)\prod\limits_{\substack{1 \leq j \leq N \\ j \neq i }}(1-\frac{y_j}{y_i})} .
\end{multline}

For our purpose, we will make use of the following particular form of \eqref{eqab3}, which we arrive by letting $x_i \to \frac{1}{x_i} , 1 \leq i \leq M$ and $y_i \to \frac{1}{y_i}, 1 \leq i \leq N$.

\begin{multline}
    \tag{1.5} \label{eqab4}
 \frac{1}{\prod\limits_{1 \leq i \leq M}(x_i - t)\prod\limits_{1 \leq i \leq N}(y_i-\frac{1}{t})} = \sum\limits_{1 \leq i \leq M}\frac{1}{(x_i - t)\prod\limits_{1 \leq j \leq N}(y_j -\frac{1}{x_i})\prod\limits_{\substack{1 \leq j \leq M \\ j \neq i }}(x_j-x_i)} \\ +\sum\limits_{1 \leq i \leq N}\frac{1}{y_i(y_i t-1)\prod\limits_{1 \leq j \leq M}(x_j- \frac{1}{y_i})\prod\limits_{\substack{1 \leq j \leq N \\ j \neq i }}(y_j-y_i)} ,
\end{multline}
which holds provided that each $x_i$, $1 \leq i \leq M$ are distinct and each $y_i$, $1 \leq i \leq N$ are distinct.

And let us recall the following two Mittag-Leffler expansion formulae :
\[ \tag{1.6} \label{eqb2}
\frac{\pi \cos(\theta x)}{2x\sin(\pi x)} - \frac{1}{2x^2}= \sum\limits_{k=1}^{\infty}\frac{(-1)^{k-1}\cos(k \theta)}{k^2-x^2},\hspace{0.12cm} |\theta| \leq \pi,
\]

provided $x \neq \pm r$, where $r$ is a natural number.

\[ \tag{1.7} \label{eqb3}
\frac{\pi\sin(\theta x)}{4x\cos(\frac{\pi x}{2})} = \sum\limits_{k=0}^{\infty} \frac{(-1)^k \sin((2k+1)\theta)}{(2k+1)^2-x^2} , \hspace{0.12cm} |\theta| \leq \frac{\pi}{2},
\]

provided $x \neq \pm r$, where $r$ is odd. These two expansions also appear in Ramanujan's notebook (see [1, p. 237]). For our purpose the following forms will prove to be more useful :

\[ \tag{1.8} \label{eqb4}
\frac{\pi \cos(\frac{\theta x}{y})}{2xy\sin(\frac{\pi x}{y})} - \frac{1}{2x^2}= \sum\limits_{k=1}^{\infty}\frac{(-1)^{k-1}\cos(k \theta)}{k^2y^2-x^2},\hspace{0.12cm} |\theta| \leq \pi,
\]

\[ \tag{1.9} \label{eqb5}
\frac{\pi\sin(\frac{\theta x}{y})}{4xy\cos(\frac{\pi x}{2y})} = \sum\limits_{k=0}^{\infty} \frac{(-1)^k \sin((2k+1)\theta)}{(2k+1)^2y^2-x^2} , \hspace{0.12cm}  |\theta| \leq \frac{\pi}{2}.
\]

Our main objective in this paper is to prove the following four propositions utilising the partial fraction identities aforementioned.
\vspace{0.7cm}

\begin{theorem}
For $|\theta_i| \leq \pi$, $1 \leq i \leq M$ and arbitrary $t$, and if for complex values $a_1,...,a_M$ and nonzero complex values $x_1,...,x_M$ such that for all $i \neq j$, and for all natural numbers $n$ and $k$, we have $x_i^2n^2-x_j^2k^2 \neq a_j^2-a_i^2$, and if for each $i \neq j$, $\textbf{Im}(\frac{x_i}{x_j}) \neq 0$ or $x_i = x_j$,  there holds :
\begin{multline} \tag{1.10} \label{eqt1} \sum\limits_{1 \leq i \leq M} \hspace{0.12cm}\sum\limits_{n=1}^{\infty}\frac{(-1)^{n-1}\cos(n \theta_i)}{x_i^2n^2+a_i^2-t^2}\prod\limits_{\substack{1\leq j \leq M \\ j \neq i}}\left\{\frac{\pi }{2x_j\sqrt{x_i^2 n^2 + a_i^2-a_j^2}}\frac{\cos{\left(\frac{\theta_j \sqrt{x_i^2 n^2+a_i^2-a_j^2}}{x_j} \right)}}{\sin{\left(\frac{\pi \sqrt{x_i^2 n^2+a_i^2-a_j^2}}{x_j} \right)}}- \frac{1}{2(x_i^2 n^2+a_i^2-a_j^2)} \right\}  \\ =\prod\limits_{1\leq i \leq M}\left\{\frac{\pi }{2x_i\sqrt{t^2-a_i^2}}\frac{\cos{\left(\frac{\theta_i \sqrt{t^2-a_i^2}}{x_i} \right)}}{\sin{\left(\frac{\pi \sqrt{t^2-a_i^2}}{x_i} \right)}} - \frac{1}{2(t^2-a_i^2)}\right\}. 
\end{multline}
\end{theorem}
\vspace{0.7cm}
\begin{theorem}
For $|\theta_i| \leq \frac{\pi}{2}$, $1 \leq i \leq M$ and arbitrary $t$, and for complex values $a_1,...,a_M$ and nonzero complex values $x_1,...,x_M$ such that for all $i \neq j$, and for all odd natural numbers $n$ and $k$, we have $x_i^2n^2-x_j^2k^2 \neq a_j^2-a_i^2$, and if for each $i \neq j$, $\textbf{Im}(\frac{x_i}{x_j}) \neq 0$ or $x_i = x_j$ , there holds :
\begin{multline}\tag{1.11} \label{eqt2}
\sum\limits_{1 \leq i \leq M} x_i\sum\limits_{n=0}^{\infty}
\frac{(-1)^n\sin((2n+1)\theta_i)}{x_i^2(2n+1)^2+a_i^2-t^2}\prod\limits_{\substack{1 \leq j \leq M \\ j \neq i }}\frac{\sin\left(\frac{\theta_j \sqrt{x_i^2(2n+1)^2+a_i^2-a_j^2}}{x_j}\right)}{\sqrt{x_i^2(2n+1)^2+a_i^2-a_j^2}\cos\left(\frac{\pi \sqrt{x_i^2(2n+1)^2+a_i^2-a_j^2}}{2x_j}\right)} \\= \frac{\pi}{4}\prod\limits_{1\leq i \leq M}\frac{\sin\left(\frac{\theta_i \sqrt{t^2-a_i^2}}{x_i}\right)}{\sqrt{t^2-a_i^2}\cos\left(\frac{\pi \sqrt{t^2-a_i^2}}{2x_i}\right)}.
\end{multline}
\end{theorem}
\vspace{0.7cm}
\begin{theorem}
For $|\theta_i| \leq \frac{\pi}{2}$, $1 \leq i \leq M$, $|\beta_i| \leq \frac{\pi}{2}$, $1 \leq i \leq N$ and arbitrary $t$, and for complex values $a_1,...,a_M$ and $b_1,...,b_N$, nonzero complex values $x_1,...,x_M$ and $y_1,...,y_N$ such that for all odd natural numbers $n$ and $k$, any two distinct $x_i$ and $x_j$, $a_i$ and $a_j$, and for any two distinct $y_i$ and $y_j$, $b_i$ and $b_j$, we have $x_i^2n^2-x_j^2k^2 \neq a_j^2-a_i^2$ and $y_i^2n^2-y_j^2k^2 \neq b_j^2-b_i^2$, and if for each $1 \leq i,j \leq M$, $i \neq j$, we have $\textbf{Im}(\frac{x_i}{x_j} \neq 0)$ or $x_i = x_j$, and for each $1 \leq i,j \leq N$, $i \neq j$, we have $\textbf{Im}(\frac{y_i}{y_j} \neq 0)$ or $y_i = y_j$,there holds :  

\begin{multline}
    \sum\limits_{1 \leq i \leq M} x_i \sum\limits_{n=0}^{\infty}\frac{(-1)^{n}\sin((2n+1)\theta_i)}{(2n+1)^2x_i^2+a_i^2-t^2} \prod\limits_{1 \leq j \leq N}\frac{\sinh\left(\frac{\beta_j \sqrt{b_j^2 -\frac{1}{(2n+1)^2x_i^2+a_i^2}}}{y_j}\right)}{\sqrt{b_j^2 -\frac{1}{(2n+1)^2x_i^2+a_i^2}}\cosh\left(\frac{\pi \sqrt{b_j^2 -\frac{1}{(2n+1)^2x_i^2+a_i^2}}}{2y_j}\right)}\prod\limits_{\substack{1 \leq j \leq M \\ j \neq i }}\frac{\sin(\frac{\theta_j \sqrt{(2n+1)^2x_i^2+a_i^2-a_j^2}}{x_j})}{\sqrt{(2n+1)^2x_i^2+a_i^2-a_j^2}\cos(\frac{\pi \sqrt{(2n+1)^2x_i^2+a_i^2-a_j^2}}{2x_j})} \\+ \sum\limits_{1 \leq i \leq N} y_i \sum\limits_{n=0}^{\infty}\frac{(-1)^{n}\sin((2n+1)\beta_i)}{((2n+1)^2y_i^2+b_i^2)((2n+1)^2y_i^2t^2+b_i^2t^2-1)} \prod\limits_{1 \leq j \leq M}\frac{\sinh\left(\frac{\theta_j \sqrt{a_j^2 -\frac{1}{(2n+1)^2y_i^2+b_i^2}}}{x_j}\right)}{\sqrt{a_j^2 -\frac{1}{(2n+1)^2y_i^2+b_i^2}}\cosh\left(\frac{\pi \sqrt{a_j^2 -\frac{1}{(2n+1)^2y_i^2+b_i^2}}}{2x_j}\right)} \times\\ \tag{1.12} \label{eqab8} \times \prod\limits_{\substack{1 \leq j \leq N \\ j \neq i }}\frac{\sin(\frac{\beta_j \sqrt{(2n+1)^2y_i^2+b_i^2-b_j^2}}{y_j})}{\sqrt{(2n+1)^2y_i^2+b_i^2-b_j^2}\cos(\frac{\pi \sqrt{(2n+1)^2y_i^2+b_i^2-b_j^2}}{2y_j})} = \frac{\pi}{4}\prod\limits_{1 \leq i \leq M}\frac{\sin(\frac{\theta_i \sqrt{t^2-a_i^2}}{x_i})}{\sqrt{t^2-a_i^2}\cos(\frac{\pi \sqrt{t^2-a_i^2}}{2x_i})} \prod\limits_{1 \leq i \leq N} \frac{\sinh(\frac{\beta_i \sqrt{b_i^2-\frac{1}{t^2}}}{y_i})}{\sqrt{b_i^2-\frac{1}{t^2}}\cosh(\frac{\pi \sqrt{b_i^2-\frac{1}{t^2}}}{2y_i})}.
\end{multline}
\end{theorem}

It will be evident that \textbf{Proposition 1} entails the famous formula of Ramanujan for odd zeta values and that \textbf{Proposition 2} entails another identity (see (18.2) in [1, p.269]) by Ramanujan. We shall also provide a general procedure for determining identities of this type. In Section 3, investigating on a special case of \textbf{Proposition 1}, we shall prove the following generalization of the transformation formula for the Dedekind eta function.

\vspace{0.7cm}
\begin{theorem}
For arbitrary $w$ and for complex values $x_1$, $x_2$, with $\textbf{Re}(x_1), \textbf{Re}(x_2) >0$, satisfying $x_2^2n^2+x_1^2k^2 \neq -w^2$ for all natural numbers $n$ and $k$, there holds :
\begin{multline}
    \frac{(1-e^{-\frac{2\pi}{x_2}\sqrt{x_1^2+w^2}})(1-e^{-\frac{2\pi}{x_2}\sqrt{4x_1^2+w^2}})(1-e^{-\frac{2\pi}{x_2}\sqrt{9x_1^2+w^2}}) \hspace{0.07cm}.\hspace{0.05cm}. \hspace{0.05cm}.}{(1-e^{-\frac{2\pi}{x_1}\sqrt{x_2^2+w^2}})(1-e^{-\frac{2\pi}{x_1}\sqrt{4x_2^2+w^2}})(1-e^{-\frac{2\pi}{x_1}\sqrt{9x_2^2+w^2}})\hspace{0.07cm}.\hspace{0.05cm}. \hspace{0.05cm}. } \\ =  \left(\frac{x_1}{x_2}\right)^{ \frac{\pi w^2}{ 2 x_1 x_2}}\sqrt{\frac{\sinh(\frac{\pi w}{x_1})}{\sinh(\frac{\pi w}{x_2})}} e^{\frac{\pi}{12}(\frac{x_1}{x_2}-\frac{x_2}{x_1})+ \pi \sum\limits_{n=1}^{\infty} \left\{ \frac{\sqrt{x_2^2n^2+w^2}}{x_1}-\frac{\sqrt{x_1^2n^2+w^2}}{x_2} - \frac{x_2n}{x_1}+ \frac{x_1n}{x_2}\right\}} . \tag{1.13}
\end{multline}
\end{theorem}

When $w = 0$, this reduces to the transformation formula for the Dedekind eta function. In Section 4, we shall unearth some further consequences of our results. In Section 5, we shall derive a formula relating to a general series akin to the generalized Lambert series studied by Dixit and Maji [6], and Kanemitsu, Tanigawa and Yoshimoto [7].

\section{Derivations}

We will first try to prove \textbf{Proposition 1}. Let us make the transformations $x_i \to x_i^2n_i^2+a_i^2$ and $t \to t^2$ in \eqref{eqb1}, with each $x_i \neq 0$, to get :

\[ \tag{2.1} \label{eqb6}\frac{1}{\prod\limits_{1\leq i \leq M}(x_i^2n_i^2-(t^2-a_i^2))}  =\sum\limits_{1 \leq i \leq M} \frac{1}{(x_i^2n_i^2+a_i^2-t^2)\prod\limits_{\substack{1 \leq j \leq M \\ j \neq i }}(x_j^2n_j^2-(x_i^2n_i^2+a_i^2-a_j^2))}.\]

Multiply both sides by $\prod\limits_{i=1}^{M} (-1)^{n_i-1}\cos(n_i\theta_i)$, where $|\theta_i| \leq \pi$, $1 \leq i \leq M$, to get :

\[ \tag{2.2} \label{eqb7} \prod\limits_{1\leq i \leq M}\frac{(-1)^{n_i-1}\cos(n_i\theta_i)}{x_i^2n_i^2-(t^2-a_i^2)}  =\sum\limits_{1 \leq i \leq M} \frac{(-1)^{n_i-1}\cos(n_i\theta_i)}{x_i^2n_i^2+a_i^2-t^2} \prod\limits_{\substack{1 \leq j \leq M \\ j \neq i }}\frac{(-1)^{n_j-1}\cos(n_j\theta_j)}{x_j^2n_j^2-(x_i^2n_i^2+a_i^2-a_j^2)}.\]

Now let us sum both sides of this equality for each $n_k, 1 \leq k \leq M$ over all natural numbers so that the sum over the $'n_k \hspace{0.025cm}'$s is $\sum\limits_{n_1=1}^{\infty}...\sum\limits_{n_M=1}^{\infty}$, then with \eqref{eqb4} in mind the left hand side will resolve into :

\[\prod\limits_{1\leq i \leq M}\left\{\frac{\pi }{2x_i\sqrt{t^2-a_i^2}}\frac{\cos{\left(\frac{\theta_i \sqrt{t^2-a_i^2}}{x_i} \right)}}{\sin{\left(\frac{\pi \sqrt{t^2-a_i^2}}{x_i} \right)}} - \frac{1}{2(t^2-a_i^2)}\right\}.\]

For the right hand side, with respect to each term $\frac{(-1)^{n_i-1}\cos(n_i\theta_i)}{x_i^2n_i^2+a_i^2-t^2} \prod\limits_{\substack{1 \leq j \leq M \\ j \neq i }}\frac{(-1)^{n_j-1}\cos(n_j\theta_j)}{(x_j^2n_j^2-(x_i^2n_i^2+a_i^2-a_j^2))}$, let us rearrange the order of the corresponding summation over $n_k, 1 \leq k \leq M$ as $\sum\limits_{n_i=1}^{\infty}\sum\limits_{n_1=1}^{\infty}...\sum\limits_{n_{i-1}=1}^{\infty}\sum\limits_{n_{i+1}=1}^{\infty}...\sum\limits_{n_M=1}^{\infty}$ so that the sum over $n_i$ will be the outermost sum. Then, again with \eqref{eqb4} in mind, the right hand side resolves into :

\[\sum\limits_{1 \leq i \leq M} \hspace{0.12cm}\sum\limits_{n=1}^{\infty}\frac{(-1)^{n-1}\cos(n \theta_i)}{x_i^2n^2+a_i^2-t^2}\prod\limits_{\substack{1\leq j \leq M \\ j \neq i}}\left\{\frac{\pi }{2x_j\sqrt{x_i^2 n^2 + a_i^2-a_j^2}}\frac{\cos{\left(\frac{\theta_j \sqrt{x_i^2 n^2+a_i^2-a_j^2}}{x_j} \right)}}{\sin{\left(\frac{\pi \sqrt{x_i^2 n^2+a_i^2-a_j^2}}{x_j} \right)}}- \frac{1}{2(x_i^2 n^2+a_i^2-a_j^2)} \right\}.\] 

Thus, we arrived at \textbf{Proposition 1}. The decomposition in \eqref{eqb6} is justified if $x_i^2n^2-x_j^2k^2 \neq a_j^2-a_i^2$ for all $i\neq j$. Rearrangements of the order of summation are justified if each of the series 

\[
\sum\limits_{n=1}^{\infty}\frac{(-1)^{n-1}\cos(n \theta_i)}{x_i^2n^2+a_i^2-t^2}\prod\limits_{\substack{1\leq j \leq M \\ j \neq i}}\left\{\frac{\pi }{2x_j\sqrt{x_i^2 n^2 + a_i^2-a_j^2}}\frac{\cos{\left(\frac{\theta_j \sqrt{x_i^2 n^2+a_i^2-a_j^2}}{x_j} \right)}}{\sin{\left(\frac{\pi \sqrt{x_i^2 n^2+a_i^2-a_j^2}}{x_j} \right)}}- \frac{1}{2(x_i^2 n^2+a_i^2-a_j^2)} \right\},
\]

converges absolutely for all $i$. Each of the above series converges absolutely for all $i$ when the sequence

\[
\left\{\frac{1}{\sqrt{x_i^2 n^2 + a_i^2-a_j^2}\sin{\left(\frac{\pi \sqrt{x_i^2 n^2+a_i^2-a_j^2}}{x_j} \right)}} \right\}_{n \in \mathbb{N}},
\]

tends to a finite limit for all $i$ and $j$ such that $i \neq j$. And the above sequence will have a finite limit for all $i \neq j$ only if for each $i \neq j$, we have $x_i = x_j$ or the imaginary part of $\frac{x_i}{x_j}$ is nonzero ($\textbf{\textit{Im}}(\frac{x_i}{x_j}) \neq 0$).

Now we proceed to prove \textbf{Proposition 2} in the same way. Make the transformations $x_i \to x_i^2(2n_i+1)^2+a_i^2$ and $t \to t^2$ in \eqref{eqb1}, with each $x_i \neq 0$, and then  multiply both sides of the equality by $\prod\limits_{i=1}^{M} (-1)^{n_i}\sin((2n_i+1)\theta_i)$, where $|\theta_i| \leq \frac{\pi}{2}$, $1 \leq i \leq M$, to get :
\[ \tag{2.3} \label{eqb8} \prod\limits_{1\leq i \leq M}\frac{(-1)^{n_i}\sin((2n_i+1)\theta_i)}{x_i^2(2n_i+1)^2-(t^2-a_i^2)}  =\sum\limits_{1 \leq i \leq M} \frac{(-1)^{n_i}\sin((2n_i+1)\theta_i)}{x_i^2(2n_i+1)^2+a_i^2-t^2} \prod\limits_{\substack{1 \leq j \leq M \\ j \neq i }}\frac{(-1)^{n_j}\sin((2n_j+1)\theta_j)}{(x_j^2(2n_j+1)^2-(x_i^2(2n_i+1)^2+a_i^2-a_j^2))}.\]

Then, let us sum both sides of the equality for each $n_k, 1 \leq k \leq M$ over all non-negative integers $0 \leq n_k \leq \infty$. Now proceeding in the same manner as in deriving \textbf{Proposition 1}, this time with \eqref{eqb5} in mind, we arrive at \textbf{Proposition 2}. The manipulations are justified when for all $i \neq j$, and for all odd natural numbers $n$ and $k$, we have $x_i^2n^2-x_j^2k^2 \neq a_j^2-a_i^2$, and if for each $i \neq j$, $\textbf{\textit{Im}}(\frac{x_i}{x_j}) \neq 0$ or $x_i = x_j$. \\
\vspace{0.7cm}
From these two derivations, we see that this general method or procedure can be done on any collection of $M$ functions $\{f_i\}_{1 \leq i \leq M}$ such that if we let each $f_i$ to have an expansion of the form :

\[ \tag{2.4} \label{eqb9}
f_i(z) = \sum\limits_{n=1}^{\infty}\frac{R_{i}(n)}{c_{i}(n)-z}.
\]
 And since by \eqref{eqb1}, it holds that : 
\[ \tag{2.5} \label{eqb10} \prod\limits_{1\leq i \leq M}\frac{R_i(n_i)}{x_ic_{i}(n_i)-(t-a_i)}  =\sum\limits_{1 \leq i \leq M} \frac{R_i(n_i)}{(x_ic_i(n_i)+a_i-t)}\prod\limits_{\substack{1 \leq j \leq M \\ j \neq i }}\frac{R_j(n_j)}{x_jc_j(n_j)-(x_ic_i(n_i)+a_i-a_j)}.\]
Using the same form of argument as the cases before, we arrive to the following generalization.
\vspace{0.7cm}

\begin{theorem}
    Let $\{f_i\}_{1 \leq i \leq M}$ be a collection of $M$ functions with each $f_i$ having an expansion of the form :
\[
f_i(z) = \sum\limits_{n=1}^{\infty}\frac{R_{i}(n)}{c_{i}(n)-z}.
\]
Then, for arbitrary $t$ and for complex values $x_1,...,x_M$ and $a_1,...,a_M$ such that for any two distinct $x_i$ and $x_j$, $a_i$ and $a_j$,  $x_ic_i(n)-x_jc_j(k) \neq a_j-a_i$ for all natural numbers $n$ and $k$, there holds :
    \[\tag{2.6} \label{eqt3}
        \prod\limits_{1\leq i \leq M}f_i\left(\frac{t-a_i}{x_i}\right) = \sum\limits_{1 \leq i \leq M} x_i\sum\limits_{n=1}^{\infty}\frac{R_i(n)}{x_ic_i(n)+a_i-t}\prod\limits_{\substack{1 \leq j \leq M \\ j \neq i }} f_j\left(\frac{x_ic_i(n)+a_i-a_j}{x_j}\right).
    \]
    Provided that each of the series from the right hand side converges absolutely. 
\end{theorem}

Now we proceed to derive \textbf{Proposition 3}. In \eqref{eqab4} make the transformations $t \to t^2$, \hspace{0.0025cm} $x_i \to (2m_i+1)^2x_i^2+a_i^2 , 1 \leq i \leq M$ and $y_i \to (2n_i+1)^2y_i^2+b_i^2, 1 \leq i \leq N$, with each $x_i \neq 0$ and $y_i \neq 0$, then we multiply both sides of the equality by $\prod\limits_{1 \leq i \leq M} (-1)^{m_i}\sin((2m_i+1)\theta_i) \prod\limits_{1 \leq i \leq N} (-1)^{n_i}\sin((2n_i+1)\beta_i)$, where $|\theta_i| \leq \frac{\pi}{2}$, $1 \leq i \leq M$, $|\beta_i| \leq \frac{\pi}{2}$, $1 \leq i \leq N$, to get : 

\begin{multline}
\prod\limits_{1 \leq i \leq M} \frac{(-1)^{m_i}\sin((2m_i+1)\theta_i)}{(2m_i+1)^2 x_i^2+a_i^2-t^2} \prod\limits_{1 \leq i \leq N} \frac{(-1)^{n_i}\sin((2n_i+1)\beta_i)}{(2n_i+1)^2y_i^2+b_i^2-\frac{1}{t^2}} \\= \sum\limits_{1 \leq i \leq M}\frac{(-1)^{m_i}\sin((2m_i+1)\theta_i)}{(2m_i+1)^2x_i^2+a_i^2-t^2} \prod\limits_{1 \leq j \leq N} \frac{(-1)^{n_j}\sin((2n_j+1)\beta_j)}{(2n_j+1)^2y_j^2+b_j^2 -\frac{1}{(2m_i+1)^2x_i^2+a_i^2}}\prod\limits_{\substack{1 \leq j \leq M \\ j \neq i }}\frac{(-1)^{m_j}\sin((2m_j+1)\theta_j)}{(2m_j+1)^2x_j^2-((2m_i+1)^2x_i^2+a_i^2-a_j^2)}  \\ +\sum\limits_{1 \leq i \leq N}\frac{(-1)^{n_i}\sin((2n_i+1)\beta_i)}{((2n_i+1)^2y_i^2+b_i^2)((2n_i+1)^2y_i^2t^2+b_i^2t^2-1)} \prod\limits_{1 \leq j \leq M} \frac{(-1)^{m_j}\sin((2m_j+1)\theta_j)}{(2m_j+1)^2x_j^2+a_j^2 -\frac{1}{(2n_i+1)^2y_i^2+b_i^2}}\prod\limits_{\substack{1 \leq j \leq N \\ j \neq i }}\frac{(-1)^{n_j}\sin((2n_j+1)\beta_j)}{(2n_j+1)^2y_j^2-((2n_i+1)^2y_i^2+b_i^2-b_j^2)} .\\ \tag{2.7} \label{eqab5}
\end{multline}

Now let us sum both sides of this equality for each $m_k, 1\leq k \leq M$ and for each $n_h, 1 \leq h \leq N$ over all non-negative integers. So that the sum over the $'m_k\hspace{0.025cm}'$s and $'n_h\hspace{0.025cm}'$s is $\sum\limits_{m_1=0}^{\infty}...\sum\limits_{m_M=0}^{\infty}\sum\limits_{n_1=0}^{\infty}...\sum\limits_{n_N=0}^{\infty}$. Then with \eqref{eqb5} in mind the left hand side resolves into :
\[
\prod\limits_{1 \leq i \leq M}\frac{\pi\sin(\frac{\theta_i \sqrt{t^2-a_i^2}}{x_i})}{4x_i\sqrt{t^2-a_i^2}\cos(\frac{\pi \sqrt{t^2-a_i^2}}{2x_i})} \prod\limits_{1 \leq i \leq N} \frac{\pi\sinh(\frac{\beta_i \sqrt{b_i^2-\frac{1}{t^2}}}{y_i})}{4y_i\sqrt{b_i^2-\frac{1}{t^2}}\cosh(\frac{\pi \sqrt{b_i^2-\frac{1}{t^2}}}{2y_i})}.
\]

For the right hand side let us look at each term inside the first summation :

\[
\frac{(-1)^{m_i}\sin((2m_i+1)\theta_i)}{(2m_i+1)^2x_i^2+a_i^2-t^2} \prod\limits_{1 \leq j \leq N} \frac{(-1)^{n_j}\sin((2n_j+1)\beta_j)}{(2n_j+1)^2y_j^2+b_j^2 -\frac{1}{(2m_i+1)^2x_i^2+a_i^2}}\prod\limits_{\substack{1 \leq j \leq M \\ j \neq i }}\frac{(-1)^{m_j}\sin((2m_j+1)\theta_j)}{(2m_j+1)^2x_j^2-((2m_i+1)^2x_i^2+a_i^2-a_j^2)},
\]
and for each of these terms respectively we shall rearrange the order of summation over the $'m_k\hspace{0.025cm}'$s and $'n_h\hspace{0.025cm}'$s as $\sum\limits_{m_i=0}^{\infty}\sum\limits_{m_1=0}^{\infty}...\sum\limits_{m_{i-1}=0}^{\infty}\sum\limits_{m_{i+1}=0}^{\infty}...\sum\limits_{m_M=0}^{\infty}\sum\limits_{n_1=0}^{\infty}...\sum\limits_{n_N=0}^{\infty}$ so that the sum over $m_i$ will be the outermost sum. Then, in view of \eqref{eqb5}, we see that the first summation from the right hand side resolves into :

\begin{multline} \notag
    \sum\limits_{1 \leq i \leq M} \sum\limits_{n=0}^{\infty}\frac{(-1)^{n}\sin((2n+1)\theta_i)}{(2n+1)^2x_i^2+a_i^2-t^2} \prod\limits_{1 \leq j \leq N}\frac{\pi\sinh\left(\frac{\beta_j \sqrt{b_j^2 -\frac{1}{(2n+1)^2x_i^2+a_i^2}}}{y_j}\right)}{4y_j\sqrt{b_j^2 -\frac{1}{(2n+1)^2x_i^2+a_i^2}}\cosh\left(\frac{\pi \sqrt{b_j^2 -\frac{1}{(2n+1)^2x_i^2+a_i^2}}}{2y_j}\right)}\prod\limits_{\substack{1 \leq j \leq M \\ j \neq i }}\frac{\pi\sin(\frac{\theta_j \sqrt{(2n+1)^2x_i^2+a_i^2-a_j^2}}{x_j})}{4x_j\sqrt{(2n+1)^2x_i^2+a_i^2-a_j^2}\cos(\frac{\pi \sqrt{(2n+1)^2x_i^2+a_i^2-a_j^2}}{2x_j})}.
\end{multline}

In the same manner, we see that the second summation from the right hand side resolves into :

\begin{multline} \notag
    \sum\limits_{1 \leq i \leq N} \sum\limits_{n=0}^{\infty}\frac{(-1)^{n}\sin((2n+1)\beta_i)}{((2n+1)^2y_i^2+b_i^2)((2n+1)^2y_i^2t^2+b_i^2t^2-1)} \prod\limits_{1 \leq j \leq M}\frac{\pi\sinh\left(\frac{\theta_j \sqrt{a_j^2 -\frac{1}{(2n+1)^2y_i^2+b_i^2}}}{x_j}\right)}{4x_j\sqrt{a_j^2 -\frac{1}{(2n+1)^2y_i^2+b_i^2}}\cosh\left(\frac{\pi \sqrt{a_j^2 -\frac{1}{(2n+1)^2y_i^2+b_i^2}}}{2x_j}\right)} \times\\ \times \prod\limits_{\substack{1 \leq j \leq N \\ j \neq i }}\frac{\pi\sin(\frac{\beta_j \sqrt{(2n+1)^2y_i^2+b_i^2-b_j^2}}{y_j})}{4y_j\sqrt{(2n+1)^2y_i^2+b_i^2-b_j^2}\cos(\frac{\pi \sqrt{(2n+1)^2y_i^2+b_i^2-b_j^2}}{2y_j})}.
\end{multline}

Thus, we arrive at \textbf{Proposition 3}. The decomposition in \eqref{eqab5} is justified if for all odd $n$ and $k$, we have $x_i^2n^2-x_j^2k^2 \neq a_j^2-a_i^2$ for $1\leq i,j \leq M$, $i \neq j$ and $y_i^2n^2-y_j^2k^2 \neq b_j^2-b_i^2$ for $1\leq i,j \leq N$, $i \neq j$. The rearrangements of the order of summation are justified when each of the series 

\[
\sum\limits_{n=0}^{\infty}\frac{(-1)^{n}\sin((2n+1)\theta_i)}{(2n+1)^2x_i^2+a_i^2-t^2} \prod\limits_{1 \leq j \leq N}\frac{\pi\sinh\left(\frac{\beta_j \sqrt{b_j^2 -\frac{1}{(2n+1)^2x_i^2+a_i^2}}}{y_j}\right)}{4y_j\sqrt{b_j^2 -\frac{1}{(2n+1)^2x_i^2+a_i^2}}\cosh\left(\frac{\pi \sqrt{b_j^2 -\frac{1}{(2n+1)^2x_i^2+a_i^2}}}{2y_j}\right)}\prod\limits_{\substack{1 \leq j \leq M \\ j \neq i }}\frac{\pi\sin(\frac{\theta_j \sqrt{(2n+1)^2x_i^2+a_i^2-a_j^2}}{x_j})}{4x_j\sqrt{(2n+1)^2x_i^2+a_i^2-a_j^2}\cos(\frac{\pi \sqrt{(2n+1)^2x_i^2+a_i^2-a_j^2}}{2x_j})},
\]

\begin{multline*}
\sum\limits_{n=0}^{\infty}\frac{(-1)^{n}\sin((2n+1)\beta_i)}{((2n+1)^2y_i^2+b_i^2)((2n+1)^2y_i^2t^2+b_i^2t^2-1)} \prod\limits_{1 \leq j \leq M}\frac{\pi\sinh\left(\frac{\theta_j \sqrt{a_j^2 -\frac{1}{(2n+1)^2y_i^2+b_i^2}}}{x_j}\right)}{4x_j\sqrt{a_j^2 -\frac{1}{(2n+1)^2y_i^2+b_i^2}}\cosh\left(\frac{\pi \sqrt{a_j^2 -\frac{1}{(2n+1)^2y_i^2+b_i^2}}}{2x_j}\right)} \times\\ \times \prod\limits_{\substack{1 \leq j \leq N \\ j \neq i }}\frac{\pi\sin(\frac{\beta_j \sqrt{(2n+1)^2y_i^2+b_i^2-b_j^2}}{y_j})}{4y_j\sqrt{(2n+1)^2y_i^2+b_i^2-b_j^2}\cos(\frac{\pi \sqrt{(2n+1)^2y_i^2+b_i^2-b_j^2}}{2y_j})},
\end{multline*}
converges absolutely. They converge absolutely when for each $1 \leq i,j \leq M$, $i \neq j$, we have $\textbf{\textit{Im}}(\frac{x_i}{x_j}) \neq 0$ or $x_i = x_j$ and for each $1 \leq i,j \leq N$, $i \neq j$, we have $\textbf{\textit{Im}}(\frac{y_i}{y_j}) \neq 0$ or $y_i = y_j$. \\ 
\vspace{0.7cm}
We again note that this procedure for deriving \textbf{Proposition 3} can be done on any two collections of $M$ functions $\{f_i\}_{1 \leq i \leq M}$ and $N$ functions $\{g_i\}_{1 \leq i \leq N}$ such that if each $f_k$ and $g_h$ have expansions of the form :

\[ \tag{2.8} \label{eqab6}
f_k(z) = \sum\limits_{n=1}^{\infty}\frac{R_{k}(n)}{c_{k}(n)-z},
\]

\[ \tag{2.9} \label{eqab7}
g_h(z) = \sum\limits_{n=1}^{\infty}\frac{P_{h}(n)}{l_{h}(n)-z}.
\]
And since by \eqref{eqab4}, it holds that : 

\begin{multline}
    \tag{2.10} \label{eqab9}
 \prod\limits_{1 \leq i \leq M}\frac{R_i(m_i)}{x_ic_i(m_i) - (t-a_i)} \prod\limits_{1 \leq i \leq N}\frac{P_i(n_i)}{y_il_i(n_i)-(\frac{1}{t}-b_i)} \\= \sum\limits_{1 \leq i \leq M}\frac{R_i(m_i)}{x_ic_i(m_i)-(t-a_i)}\prod\limits_{1 \leq j \leq N} \frac{P_j(n_j)}{y_jl_j(n_j)-(\frac{1}{x_ic_i(m_i)+a_i} - b_j)} \prod\limits_{\substack{1 \leq j \leq M \\ j \neq i }} \frac{R_j(m_j)}{x_jc_j(m_j)-(x_ic_i(m_i)+a_i-a_j)} \\ +\sum\limits_{1 \leq i \leq N} \frac{P_i(n_i)}{(y_il_i(n_i)+b_i)(y_il_i(n_i)t+b_it-1)} \prod\limits_{1 \leq j \leq M}\frac{R_j(m_j)}{x_jc_j(m_j)-(\frac{1}{y_il_i(n_i)+b_i}-a_j)} \prod\limits_{\substack{1 \leq j \leq N \\ j \neq i }} \frac{P_j(n_j)}{y_jl_j(n_j)-(y_il_i(n_i)+b_i-b_j)}.
\end{multline}

Using the same form of argument as in deriving \textbf{Proposition 3}, we arrive to the following general statement.
\begin{theorem}  Let $\{f_i\}_{1 \leq i \leq M}$ and  $\{g_i\}_{1 \leq i \leq N}$ be two collections of $M$ and $N$ functions respectively such that each $f_k$ and $g_h$ having expansions of the form :

\[ 
f_k(z) = \sum\limits_{n=1}^{\infty}\frac{R_{k}(n)}{c_{k}(n)-z},
\]

\[
g_h(z) = \sum\limits_{n=1}^{\infty}\frac{P_{h}(n)}{l_{h}(n)-z}.
\]

Then for arbitrary $t$, and for complex values $x_1,...,x_M$ and $a_1,...,a_M$,  $y_1,...,y_N$ and $b_1,...,b_N$ such that for any two distinct $x_i$ and $x_j$, $a_i$ and $a_j$, and for any two distinct $y_i$ and $y_j$, $b_i$ and $b_j$, $x_ic_i(n)-x_jc_j(k) \neq a_j-a_i$ and $y_il_i(n)-y_jl_j(k) \neq b_j-b_i$ for all natural numbers $n$ and $k$, there holds :
    \begin{multline}
        \prod\limits_{1 \leq i \leq M} f_i\left( \frac{t-a_i}{x_i}\right) \prod\limits_{1 \leq i \leq N} g_i\left( \frac{\frac{1}{t}-b_i}{y_i}\right) \\= \sum\limits_{1 \leq i \leq M} x_i\sum\limits_{n=0}^{\infty}\frac{R_i(n)}{x_ic_i(n)+a_i-t}\prod\limits_{1 \leq j \leq N} g_j\left( \frac{1}{y_j(x_ic_i(n)+a_i)} - \frac{b_j}{y_j}\right) \prod\limits_{\substack{1 \leq j \leq M \\ j \neq i }} f_j\left( \frac{x_ic_i(n)+a_i-a_j}{x_j}\right) \\ \tag{2.11}\label{eqab10}+\sum\limits_{1 \leq i \leq N} y_i \sum\limits_{n=1}^{\infty} \frac{P_i(n)}{(y_il_i(n)+b_i)(y_il_i(n)t+b_it-1)} \prod\limits_{1 \leq j \leq M} f_j\left(\frac{1}{x_j(y_il_i(n)+b_i)}-\frac{a_j}{x_j} \right) \prod\limits_{\substack{1 \leq j \leq N \\ j \neq i }} g_j \left( \frac{y_il_i(n)+b_i-b_j}{y_j}\right).
    \end{multline}
    Provided that the series from the right hand side converges absolutely.
\end{theorem}

\section{A certain transformation formula}

In this section we shall prove a general transformation formula which contains as a particular case, the following well-known transformation formula for the logarithm of the Dedekind eta function. This transformation formula is found in [5, Ch. 14, Sec. 8, Cor. (ii); Ch. 16, Entry 27(iii)] and [1, p. 256].
For $\alpha, \beta>0$ and $\alpha\beta=\pi^2$ :
\begin{align}\label{logdede} \tag{3.1}
\sum_{n=1}^{\infty}\frac{1}{n(e^{2n\alpha}-1)}-\sum_{n=1}^{\infty}\frac{1}{n(e^{2n\beta}-1)}=\frac{\beta-\alpha}{12}+\frac{1}{4}\ln\left(\frac{\alpha}{\beta}\right).
\end{align}

The following generalization is also found in [5, vol. I, p. 257, no. 12; vol. II, p. 169, no. 8 (ii)] and [1, p. 253]. A proof has been given by Lagrange [4]. For $\alpha, \beta>0$ with $\alpha\beta=\pi$, and $0 < \alpha \eta < \pi$ :

\begin{align}\label{logdede2} \tag{3.2}
2\sum_{n=1}^{\infty}\frac{\cos(2n\eta \alpha)}{n(e^{2n\alpha^2}-1)}-2\sum_{n=1}^{\infty}\frac{\cosh(2n\eta \beta)}{n(e^{2n\beta^2}-1)}= \eta^2+\frac{\beta^2-\alpha^2}{6}+\ln\left(\frac{\sin(\alpha \eta)}{\sinh(\beta \eta)}\right).
\end{align}

The transformation formula that we shall derive here via the usage of \textbf{Proposition 1}, will contain both \eqref{logdede} and \eqref{logdede2} as special cases. Let us begin by looking into the case $M=2$ of \textbf{Proposition 1}.

\begin{multline}
    \sum\limits_{n=1}^{\infty} \frac{(-1)^{n-1}\cos(n \theta_1)}{x_1^2n^2+a_1^2-t^2}\left(\frac{\pi }{2x_2\sqrt{x_1^2 n^2 + a_1^2-a_2^2}}\frac{\cos{\left(\frac{\theta_2 \sqrt{x_1^2 n^2+a_1^2-a_2^2}}{x_2} \right)}}{\sin{\left(\frac{\pi \sqrt{x_1^2 n^2+a_1^2-a_2^2}}{x_2} \right)}}- \frac{1}{2(x_1^2 n^2+a_1^2-a_2^2)} \right) \\+ \sum\limits_{n=1}^{\infty} \frac{(-1)^{n-1}\cos(n \theta_2)}{x_2^2n^2+a_2^2-t^2}\left(\frac{\pi }{2x_1\sqrt{x_2^2 n^2 + a_2^2-a_1^2}}\frac{\cos{\left(\frac{\theta_1 \sqrt{x_2^2 n^2+a_2^2-a_1^2}}{x_1} \right)}}{\sin{\left(\frac{\pi \sqrt{x_2^2 n^2+a_2^2-a_1^2}}{x_1} \right)}}- \frac{1}{2(x_2^2 n^2+a_2^2-a_1^2)} \right) \\ \tag{3.3} = \left\{\frac{\pi }{2x_1\sqrt{t^2-a_1^2}}\frac{\cos{\left(\frac{\theta_1 \sqrt{t^2-a_1^2}}{x_1} \right)}}{\sin{\left(\frac{\pi \sqrt{t^2-a_1^2}}{x_1} \right)}} - \frac{1}{2(t^2-a_1^2)}\right\}\left\{\frac{\pi }{2x_2\sqrt{t^2-a_2^2}}\frac{\cos{\left(\frac{\theta_2 \sqrt{t^2-a_2^2}}{x_2} \right)}}{\sin{\left(\frac{\pi \sqrt{t^2-a_2^2}}{x_2} \right)}} - \frac{1}{2(t^2-a_2^2)}\right\}.
\end{multline}

Let $t \to it^{-1}$, and divide both sides by $t^2$ to arrive at :
\begin{multline}
    \sum\limits_{n=1}^{\infty} \frac{(-1)^{n-1}\cos(n \theta_1)}{(x_1^2n^2+a_1^2)t^2+1}\left(\frac{\pi }{2x_2\sqrt{x_1^2 n^2 + a_1^2-a_2^2}}\frac{\cos{\left(\frac{\theta_2 \sqrt{x_1^2 n^2+a_1^2-a_2^2}}{x_2} \right)}}{\sin{\left(\frac{\pi \sqrt{x_1^2 n^2+a_1^2-a_2^2}}{x_2} \right)}}- \frac{1}{2(x_1^2 n^2+a_1^2-a_2^2)} \right) \\+ \sum\limits_{n=1}^{\infty} \frac{(-1)^{n-1}\cos(n \theta_2)}{(x_2^2n^2+a_2^2)t^2+1}\left(\frac{\pi }{2x_1\sqrt{x_2^2 n^2 + a_2^2-a_1^2}}\frac{\cos{\left(\frac{\theta_1 \sqrt{x_2^2 n^2+a_2^2-a_1^2}}{x_1} \right)}}{\sin{\left(\frac{\pi \sqrt{x_2^2 n^2+a_2^2-a_1^2}}{x_1} \right)}}- \frac{1}{2(x_2^2 n^2+a_2^2-a_1^2)} \right) \\ \tag{3.4} = \left\{-\frac{\pi }{2x_1\sqrt{1+t^2a_1^2}}\frac{\cosh{\left(\frac{\theta_1 \sqrt{1+t^2a_1^2}}{x_1t} \right)}}{\sinh{\left(\frac{\pi \sqrt{1+t^2a_1^2}}{x_1t} \right)}} + \frac{t}{2(1+a_1^2t^2)}\right\}\left\{-\frac{\pi }{2x_2\sqrt{1+t^2a_2^2}}\frac{\cosh{\left(\frac{\theta_2 \sqrt{1+t^2a_2^2}}{x_2t} \right)}}{\sinh{\left(\frac{\pi \sqrt{1+t^2a_2^2}}{x_2t} \right)}} + \frac{t}{2(1+a_2^2t^2)}\right\}.
\end{multline}
Now we shall restrict both $|\theta_1|$ and $|\theta_2|$ to be strictly less than $\pi$ and let us put $t \to 0$. Since $|\theta_1|$ and $|\theta_2|$ are strictly less than $\pi$, we see that the right hand side of the equality tends to $0$ as $t \to 0$. Thus, we arrive at :

\begin{multline}
    \sum\limits_{n=1}^{\infty} \frac{\pi (-1)^{n-1}\cos(n \theta_1)}{2x_2\sqrt{x_1^2 n^2 + a_1^2-a_2^2}}\frac{\cos{\left(\frac{\theta_2 \sqrt{x_1^2 n^2+a_1^2-a_2^2}}{x_2} \right)}}{\sin{\left(\frac{\pi \sqrt{x_1^2 n^2+a_1^2-a_2^2}}{x_2} \right)}}+\sum\limits_{n=1}^{\infty} \frac{\pi (-1)^{n-1}\cos(n \theta_2)}{2x_1\sqrt{x_2^2 n^2 + a_2^2-a_1^2}}\frac{\cos{\left(\frac{\theta_1 \sqrt{x_2^2 n^2+a_2^2-a_1^2}}{x_1} \right)}}{\sin{\left(\frac{\pi \sqrt{x_2^2 n^2+a_2^2-a_1^2}}{x_1} \right)}}\\ \tag{3.5} - \frac{1}{2}\sum\limits_{n=1}^{\infty} \frac{(-1)^{n-1}\cos(n \theta_1)}{x_1^2n^2+a_1^2-a_2^2}- \frac{1}{2}\sum\limits_{n=1}^{\infty} \frac{(-1)^{n-1}\cos(n \theta_2)}{x_2^2n^2+a_2^2-a_1^2} = 0 .
\end{multline}

Let $x_2 \to ix_2$, and let $a_1^2-a_2^2=w^2$, with \eqref{eqb4} in mind, we arrive at :

\begin{multline}
    \sum\limits_{n=1}^{\infty} \frac{ (-1)^{n-1}\cos(n \theta_1)}{x_2\sqrt{x_1^2 n^2 + w^2}}\frac{\cosh{\left(\frac{\theta_2 \sqrt{x_1^2 n^2+w^2}}{x_2} \right)}}{\sinh{\left(\frac{\pi \sqrt{x_1^2 n^2+w^2}}{x_2} \right)}}-\sum\limits_{n=1}^{\infty} \frac{ (-1)^{n-1}\cos(n \theta_2)}{x_1\sqrt{x_2^2 n^2 + w^2}}\frac{\cosh{\left(\frac{\theta_1 \sqrt{x_2^2 n^2+w^2}}{x_1} \right)}}{\sinh{\left(\frac{\pi \sqrt{x_2^2 n^2+w^2}}{x_1} \right)}} \tag{3.6} = \frac{ \cosh(\frac{\theta_2 w}{x_2})}{2wx_2 \sinh(\frac{\pi w}{x_2})} - \frac{ \cosh(\frac{\theta_1 w}{x_1})}{2wx_1 \sinh(\frac{\pi w}{x_1})}.
\end{multline}

Now, let $\theta_1 \to \pi - \theta_1$ and $\theta_2 \to \pi - \theta_2$ so that the restriction on $\theta_1$ and $\theta_2$ becomes $0 < \theta_1 < 2\pi$ and $0 < \theta_2 < 2\pi$. Now using the elementary fact $\cos(k\pi - x) = (-1)^k \cos(x)$ and the addition formula for hyperbolic cosine function, we arrive at :

\begin{multline}
    -\sum\limits_{n=1}^{\infty} \frac{ \cos(n \theta_1)}{x_2\sqrt{x_1^2 n^2 + w^2}} \left( \cosh{\left(\frac{\theta_2 \sqrt{x_1^2 n^2+w^2}}{x_2} \right)}\coth{\left(\frac{\pi \sqrt{x_1^2 n^2+w^2}}{x_2}\right)} - \sinh{\left(\frac{\theta_2 \sqrt{x_1^2 n^2+w^2}}{x_2} \right)}\right) \\ \tag{3.7}+\sum\limits_{n=1}^{\infty} \frac{ \cos(n \theta_2)}{x_1\sqrt{x_2^2 n^2 + w^2}} \left( \cosh{\left(\frac{\theta_1 \sqrt{x_2^2 n^2+w^2}}{x_1}\right)}\coth{\left(\frac{\pi \sqrt{x_2^2 n^2+w^2}}{x_1} \right)} - \sinh{\left(\frac{\theta_1 \sqrt{x_2^2 n^2+w^2}}{x_1} \right)} \right)  \\= \frac{1}{2wx_2}\cosh{(\frac{\theta_2w}{x_2})}\coth{(\frac{\pi w}{x_2})} - \frac{1}{2wx_2}\sinh{(\frac{\theta_2 w}{x_2})}  -\frac{1}{2wx_1}\cosh{(\frac{\theta_1 w}{x_1})}\coth{(\frac{\pi w}{x_1})} + \frac{1}{2wx_1}\sinh{(\frac{\theta_1 w}{x_1})} .
\end{multline}

And since, $\coth{(x)}=1+\frac{2}{e^{2x}-1}$, making some further manipulations, we arrive at :
\vspace{0.7cm}
\begin{theorem}
For $0 < \theta_1 < 2\pi$, $0 < \theta_2 < 2\pi$, and arbitrary $w$, and for complex values $x_1$, $x_2$, with $\frac{x_1}{x_2}$ not purely imaginary, satisfying $x_2^2n^2+x_1^2k^2 \neq -w^2$ for all natural numbers $n$ and $k$, there holds :
\begin{multline}
    \frac{1}{x_2} \sum\limits_{n=1}^{\infty}\frac{\cos(n 
    \theta_1) \cosh\left(\frac{\theta_2 \sqrt{x_1^2n^2+w^2}}{x_2} \right)}{\sqrt{x_1^2n^2+w^2}\left( e^{\frac{2 \pi}{x_2}\sqrt{x_1^2n^2+w^2}}-1\right)} - \frac{1}{x_1} \sum\limits_{n=1}^{\infty}\frac{\cos(n 
    \theta_2) \cosh\left(\frac{\theta_1 \sqrt{x_2^2n^2+w^2}}{x_1} \right)}{\sqrt{x_2^2n^2+w^2}\left( e^{\frac{2 \pi}{x_1}\sqrt{x_2^2n^2+w^2}}-1\right)} = \frac{\cosh(\frac{\theta_1 w}{x_1})\coth(\frac{\pi w}{x_1})}{4wx_1} - \frac{\cosh(\frac{\theta_2 w}{x_2})\coth(\frac{\pi w}{x_2})}{4wx_2} \\ \tag{3.8} + \frac{\sinh(\frac{\theta_2 w}{x_2})}{4wx_2}-\frac{\sinh(\frac{\theta_1 w}{x_1})}{4wx_1} + \frac{1}{2x_1} \sum\limits_{n=1}^{\infty} \frac{\cos(n \theta_2) e^{-\frac{\theta_1}{x_1} \sqrt{x_2^2n^2+w^2}}}{\sqrt{x_2^2n^2+w^2}} - \frac{1}{2x_2} \sum\limits_{n=1}^{\infty} \frac{\cos(n \theta_1) e^{-\frac{\theta_2}{x_2} \sqrt{x_1^2n^2+w^2}}}{\sqrt{x_1^2n^2+w^2}}.
\end{multline}
\end{theorem}

We shall see that this identity contains \eqref{logdede2} and thus also \eqref{logdede} as special cases. To see this we need to only look into the case when we let $w \to 0$. Since $\coth(x)= \frac{1}{x}+\frac{x}{3}+O(x^2)$ as $x \to 0$, we have :

\[
\lim_{w \to 0} \left\{ \frac{\cosh(\frac{\theta_1 w}{x_1})\coth(\frac{\pi w}{x_1})}{4wx_1} - \frac{\cosh(\frac{\theta_2 w}{x_2})\coth(\frac{\pi w}{x_2})}{4wx_2} \right\} = \frac{\pi}{12}(\frac{1}{x_1^2}-\frac{1}{x_2^2})+ \lim_{w \to 0} \frac{1}{4 \pi}\left( \frac{\cosh(\frac{\theta_1 w}{x_1})-\cosh(\frac{\theta_2 w}{x_2})}{w^2} \right) = \frac{\pi}{12}(\frac{1}{x_1^2}-\frac{1}{x_2^2})+ \frac{1}{8 \pi}\left( \frac{\theta_1^2}{x_1^2}-\frac{\theta_2^2}{x_2^2}\right).
\]

And since $\sum\limits_{k=1}^{\infty} \frac{\cos(k \alpha) e^{-kt}}{k} = - \frac{1}{2}\ln(1 - 2 \cos(\alpha) e^{-t}+ e^{-2t})$, we arrive at the following statement : 

 \begin{multline}
     \sum\limits_{n=1}^{\infty}\frac{\cos(n 
    \theta_1) \cosh\left(\frac{n\theta_2 x_1}{x_2} \right)}{n\left( e^{\frac{2 n \pi x_1}{x_2}}-1\right)} - \sum\limits_{n=1}^{\infty}\frac{\cos(n 
    \theta_2) \cosh\left(\frac{n\theta_1 x_2}{x_1} \right)}{n\left( e^{\frac{2 n \pi x_2}{x_1}}-1\right)} = \frac{\pi}{12}(\frac{x_2}{x_1}-\frac{x_1}{x_2})+ \frac{1}{8 \pi}\left( \frac{\theta_1^2 x_2}{x_1}-\frac{\theta_2^2 x_1}{x_2}\right) + \frac{\theta_2 x_1}{4x_2}- \frac{\theta_1 x_2}{4x_1}+\frac{1}{4} \ln\left( \frac{1-2\cos(\theta_1) e^{-\frac{\theta_2x_1}{x_2}} + e^{-\frac{2\theta_2x_1}{x_2}}  }{1-2\cos(\theta_2) e^{-\frac{\theta_1x_2}{x_1}} + e^{-\frac{2\theta_1x_2}{x_1}}  }\right). \\ \tag{3.9}
    \end{multline}
With $ \sqrt{\frac{\pi x_1}{x_2}} =\alpha $ and $\sqrt{\frac{\pi x_2}{x_1}} =\beta $, and putting $\theta_1 \to \alpha \theta_1$, $\theta_2 \to \beta \theta_2$, we state the following corollary.
\vspace{0.7cm}
\begin{corollary} For $\alpha \beta = \pi$, $ 0< \theta_1 < 2\pi$ and $0 < \theta_2 < 2\pi$, there holds :
   \begin{multline} \tag{3.10}
     \sum\limits_{n=1}^{\infty}\frac{\cos(n 
    \theta_1 \alpha) \cosh\left(n\theta_2 \alpha\right)}{n\left( e^{2 n \alpha^2}-1\right)} - \sum\limits_{n=1}^{\infty}\frac{\cos(n 
    \theta_2 \beta) \cosh\left( n\theta_1 \beta\right)}{n\left( e^{2n \beta^2}-1\right)} = \frac{\theta_1^2 - \theta_2^2}{8}+ \frac{\beta^2- \alpha^2}{12}+ \frac{1}{4}\ln\left( \frac{ \cosh(\theta_2\alpha)-\cos(\theta_1\alpha) }{\cosh(\theta_1\beta)- \cos(\theta_2 \beta)}\right).
    \end{multline}
\end{corollary}

It is evident that the limiting case $\theta_2 \to 0$ entails \eqref{logdede2}. We now proceed to prove \textbf{Proposition 4}. Intuition tells us that this should be achieved by letting both $\theta_1$ and $\theta_2$ tends to zero in $(3.8)$. The left hand side yields 

\[
\frac{1}{x_2} \sum\limits_{n=1}^{\infty}\frac{1}{\sqrt{x_1^2n^2+w^2}\left( e^{\frac{2 \pi}{x_2}\sqrt{x_1^2n^2+w^2}}-1\right)} - \frac{1}{x_1} \sum\limits_{n=1}^{\infty}\frac{1}{\sqrt{x_2^2n^2+w^2}\left( e^{\frac{2 \pi}{x_1}\sqrt{x_2^2n^2+w^2}}-1\right)},
\]

while the right hand side is 

\[
\frac{\coth(\frac{\pi w}{x_1})}{4wx_1} - \frac{\coth(\frac{\pi w}{x_2})}{4wx_2} + \lim_{\substack{\theta_1 \to 0 \\ \theta_2 \to 0}} \left\{\frac{1}{2x_1} \sum\limits_{n=1}^{\infty} \frac{\cos(n \theta_2) e^{-\frac{\theta_1}{x_1} \sqrt{x_2^2n^2+w^2}}}{\sqrt{x_2^2n^2+w^2}} - \frac{1}{2x_2} \sum\limits_{n=1}^{\infty} \frac{\cos(n \theta_1) e^{-\frac{\theta_2}{x_2} \sqrt{x_1^2n^2+w^2}}}{\sqrt{x_1^2n^2+w^2}} \right\}.
\]

By evaluating the limit $\theta_2 \to 0 $ first, the limit from the right hand side resolves into 

\[
\frac{1}{2} \sum\limits_{n=1}^{\infty} \left\{ \frac{1}{x_1\sqrt{x_2^2n^2+w^2}}- \frac{1}{x_2\sqrt{x_1^2n^2+w^2}}\right\} + \frac{1}{2 x_1 x_2} \lim_{\theta_1 \to 0} \left\{ \sum\limits_{n=1}^{\infty} \frac{e^{-\frac{\theta_1}{x_1}}\sqrt{x_2^2n^2+w^2}}{n} + \ln(2 \sin(\frac{\theta_1}{2}))\right\}.
\]

When we define $ \lim_{\theta \to 0} \left\{ \sum\limits_{n=1}^{\infty} \frac{e^{-\theta\sqrt{a^2n^2+b^2}}}{n} + \ln(2\sin(\frac{\theta}{2})) \right\} = p(a,b) $, this altogether yields 

\begin{multline}
    \frac{1}{x_2} \sum\limits_{n=1}^{\infty}\frac{1}{\sqrt{x_1^2n^2+w^2}\left( e^{\frac{2 \pi}{x_2}\sqrt{x_1^2n^2+w^2}}-1\right)} - \frac{1}{x_1} \sum\limits_{n=1}^{\infty}\frac{1}{\sqrt{x_2^2n^2+w^2}\left( e^{\frac{2 \pi}{x_1}\sqrt{x_2^2n^2+w^2}}-1\right)} \\= \frac{\coth(\frac{\pi w}{x_1})}{4wx_1} - \frac{\coth(\frac{\pi w}{x_2})}{4wx_2} + \frac{1}{2} \sum\limits_{n=1}^{\infty} \left\{ \frac{1}{x_1\sqrt{x_2^2n^2+w^2}}- \frac{1}{x_2\sqrt{x_1^2n^2+w^2}}\right\} + \frac{1}{2 x_1 x_2} p(\frac{x_2}{x_1},\frac{w}{x_1}). \notag
\end{multline}

But we see that $p(a,b)$ is independent of $b$, since 

\[
\frac{\partial p(a,b)}{\partial b} = \lim_{\theta \to 0} \left\{ \sum\limits_{n=1}^{\infty} \frac{e^{-\theta\sqrt{a^2n^2+b^2}}}{n} \frac{(-\theta b)}{\sqrt{a^2n^2+b^2}} \right\} = 0.
\]

Thus we have, $p(a,b) = p(a,0) = \lim_{\theta \to 0} \left\{ \sum\limits_{n=1}^{\infty} \frac{e^{-\theta a n}}{n} + \ln(2\sin(\frac{\theta}{2})) \right\} = \lim_{\theta \to 0} \ln\left\{ e^{\frac{\theta a}{2}}\frac{\sin(\frac{\theta}{2})}{\sinh(\frac{\theta a}{2})}\right\} = \ln(\frac{1}{a}) $. Therefore, we arrive at

\begin{multline}
    \frac{1}{x_2} \sum\limits_{n=1}^{\infty}\frac{1}{\sqrt{x_1^2n^2+w^2}\left( e^{\frac{2 \pi}{x_2}\sqrt{x_1^2n^2+w^2}}-1\right)} - \frac{1}{x_1} \sum\limits_{n=1}^{\infty}\frac{1}{\sqrt{x_2^2n^2+w^2}\left( e^{\frac{2 \pi}{x_1}\sqrt{x_2^2n^2+w^2}}-1\right)} \\= \frac{\coth(\frac{\pi w}{x_1})}{4wx_1} - \frac{\coth(\frac{\pi w}{x_2})}{4wx_2} + \frac{1}{2} \sum\limits_{n=1}^{\infty} \left\{ \frac{1}{x_1\sqrt{x_2^2n^2+w^2}}- \frac{1}{x_2\sqrt{x_1^2n^2+w^2}}\right\} + \frac{1}{2 x_1 x_2} \ln(\frac{x_1}{x_2}). \tag{3.11}
\end{multline}

Multiplying both sides by $2\pi w$ and integrating with respect to $w$, from $0$ to $w$ in both sides of the equality yields 

\begin{multline}
     \sum\limits_{n=1}^{\infty} \ln(1-e^{-\frac{2 \pi}{x_2}\sqrt{x_1^2n^2+w^2}}) - \sum\limits_{n=1}^{\infty}\ln(1-e^{-\frac{2 \pi}{x_1}\sqrt{x_2^2n^2+w^2}}) -\sum\limits_{n=1}^{\infty} \ln(1-e^{-\frac{2 \pi x_1 n}{x_2}}) + \sum\limits_{n=1}^{\infty}\ln(1-e^{-\frac{2 \pi x_2 n}{x_1}}) \\= \frac{1}{2}\ln(\sinh(\frac{\pi w}{x_1}))- \frac{1}{2}\ln(\sinh(\frac{\pi w}{x_2}))- \frac{1}{2} \ln(\frac{x_2}{x_1}) + \pi \sum\limits_{n=1}^{\infty} \left\{ \frac{\sqrt{x_2^2n^2+w^2}}{x_1}- \frac{\sqrt{x_1^2n^2+w^2}}{x_2} - \frac{x_2n}{x_1} + \frac{x_1n}{x_2}\right\} + \frac{\pi w^2}{ 2 x_1 x_2} \ln(\frac{x_1}{x_2}). \tag{3.12}
\end{multline}

Using \eqref{logdede} to evaluate the part $-\sum\limits_{n=1}^{\infty} \ln(1-e^{-\frac{2 \pi x_1 n}{x_2}}) + \sum\limits_{n=1}^{\infty}\ln(1-e^{-\frac{2 \pi x_2 n}{x_1}})$, we arrive at \textbf{Proposition 4}.

\section{Some further consequences}

Let $\theta_i = \pi$ for all $1 \leq i \leq M$ in \textbf{Proposition 1} to get the following corollary. 
\vspace{0.7cm}
\begin{corollary}
For arbitrary $t$, and if for complex values $a_1,...,a_M$ and nonzero complex values $x_1,...,x_M$ such that for all $i \neq j$, and for all natural numbers $n$ and $k$, we have $x_i^2n^2-x_j^2k^2 \neq a_j^2-a_i^2$, and if for each $i \neq j$, $\textbf{Im}(\frac{x_i}{x_j}) \neq 0$ or $x_i = x_j$,  there holds :
    \begin{multline} \tag{4.1} \label{eqt4} \sum\limits_{1 \leq i \leq M} \hspace{0.12cm}\sum\limits_{n=1}^{\infty}\frac{1}{(x_i^2n^2+a_i^2-t^2)}\prod\limits_{\substack{1\leq j \leq M \\ j \neq i}}\left\{\frac{1}{2(x_i^2 n^2+a_i^2-a_j^2)}-\frac{\pi }{2x_j\sqrt{x_i^2 n^2 + a_i^2-a_j^2}}\cot{\left(\frac{\pi \sqrt{x_i^2 n^2+a_i^2-a_j^2}}{x_j} \right)} \right\}  \\ =\prod\limits_{1\leq i \leq M}\left\{\frac{1}{2(t^2-a_i^2)}-\frac{\pi }{2x_i\sqrt{t^2-a_i^2}}\cot{\left(\frac{\pi \sqrt{t^2-a_i^2}}{x_i} \right)} \right\}. 
\end{multline}
\end{corollary}

Put $M=2$ to get: 
\begin{multline}
    \tag{4.2} \label{eqb11}
    \sum\limits_{n=1}^{\infty}\frac{1}{(x_1^2n^2+a_1^2-t^2)}\left\{\frac{1}{2(x_1^2 n^2+a_1^2-a_2^2)}-\frac{\pi }{2x_2\sqrt{x_1^2 n^2 + a_1^2-a_2^2}}\cot{\left(\frac{\pi \sqrt{x_1^2 n^2+a_1^2-a_2^2}}{x_2} \right)} \right\}\\+\sum\limits_{n=1}^{\infty}\frac{1}{(x_2^2n^2+a_2^2-t^2)}\left\{\frac{1}{2(x_2^2 n^2+a_2^2-a_1^2)}-\frac{\pi }{2x_1\sqrt{x_2^2 n^2 + a_2^2-a_1^2}}\cot{\left(\frac{\pi \sqrt{x_2^2 n^2+a_2^2-a_1^2}}{x_1} \right)} \right\} \\ = \left\{\frac{1}{2(t^2-a_1^2)}-\frac{\pi }{2x_1\sqrt{t^2-a_1^2}}\cot{\left(\frac{\pi \sqrt{t^2-a_1^2}}{x_1} \right)} \right\}\left\{\frac{1}{2(t^2-a_2^2)}-\frac{\pi }{2x_2\sqrt{t^2-a_2^2}}\cot{\left(\frac{\pi \sqrt{t^2-a_2^2}}{x_2} \right)} \right\}.
\end{multline}
Which when simplified further, is equivalent to the following statement.
\vspace{0.7cm}
\begin{corollary} For arbitrary $t$, and if for complex values $a_1,a_2$ and nonzero complex values $x_1,x_2$ such that for all natural numbers $n$ and $k$, we have $x_1^2n^2-x_2^2k^2 \neq a_2^2-a_1^2$, and if $\textbf{Im}(\frac{x_1}{x_2}) \neq 0$ or $x_1 = x_2$,  there holds :
    \begin{multline}
    \tag{4.3} \label{eqt5} 
    \frac{1}{(t^2-a_1^2)x_2\sqrt{a_1^2-a_2^2}}\cot{\left( \frac{\pi \sqrt{a_1^2-a_2^2}}{x_2}\right)}+\frac{1}{(t^2-a_2^2)x_1\sqrt{a_2^2-a_1^2}}\cot{\left( \frac{\pi \sqrt{a_2^2-a_1^2}}{x_1}\right)} \\-2\sum\limits_{n=1}^{\infty}\left\{ \frac{x_2^{-1}\cot{\left( \frac{\pi \sqrt{x_1^2 n^2 + a_1^2 - a_2^2}}{x_2}\right)}}{(x_1^2n^2+a_1^2-t^2)\sqrt{x_1^2n^2+a_1^2-a_2^2}}+\frac{x_1^{-1}\cot{\left( \frac{\pi \sqrt{x_2^2 n^2 + a_2^2 - a_1^2}}{x_1}\right)}}{(x_2^2n^2+a_2^2-t^2)\sqrt{x_2^2n^2+a_2^2-a_1^2}}\right\} \\ = \frac{\pi}{x_1x_2\sqrt{t^2-a_1^2}\sqrt{t^2-a_2^2}}\cot{\left(\frac{\pi \sqrt{t^2-a_1^2}}{x_1} \right)}\cot{\left(\frac{\pi \sqrt{t^2-a_2^2}}{x_2} \right)}.
\end{multline}
\end{corollary}
For $t=0, x_1 = x_2 =1$, this gives : 
\vspace{0.7cm}
\begin{corollary}
    For complex values $a_1,a_2$ such that $a_2^2-a_1^2$ is not integral, there holds :
    \begin{multline}
    \tag{4.4} \label{eqt6} 
    \frac{1}{a_1^2\sqrt{a_1^2-a_2^2}}\cot{\left( \pi \sqrt{a_1^2-a_2^2}\right)}+\frac{1}{a_2^2\sqrt{a_2^2-a_1^2}}\cot{\left( \pi \sqrt{a_2^2-a_1^2}\right)} \\+2\sum\limits_{n=1}^{\infty}\left\{ \frac{\cot{\left( \pi \sqrt{n^2 + a_1^2 - a_2^2}\right)}}{(n^2+a_1^2)\sqrt{n^2+a_1^2-a_2^2}}+\frac{\cot{\left( \pi \sqrt{n^2 + a_2^2 - a_1^2}\right)}}{(n^2+a_2^2)\sqrt{n^2+a_2^2-a_1^2}}\right\} \\ = -\frac{\pi}{a_1a_2}\coth{\left(\pi a_1\right)}\coth{\left(\pi a_2 \right)} .
\end{multline}
\end{corollary}

Now we shall note down the following particular cases. In \eqref{eqt6}, for $a_1 = 1, a_2 = \frac{1}{\sqrt{2}}$ :
\[ \tag{4.5} \label{eqb12}
\sum\limits_{n=1}^{\infty} \left\{ \frac{\cot(\pi\sqrt{n^2+\frac{1}{2}})}{(n^2+1)\sqrt{n^2+\frac{1}{2}}}+\frac{\cot(\pi\sqrt{n^2-\frac{1}{2}})}{(n^2+\frac{1}{2})\sqrt{n^2-\frac{1}{2}}}\right\} = \sqrt{2}\coth\left(\frac{\sqrt{2}\pi}{2}\right) - \frac{\sqrt{2}}{2}\cot\left(\frac{\sqrt{2}\pi}{2}\right)-\frac{\pi \sqrt{2}}{2}\coth(\pi)\coth\left(\frac{\sqrt{2}\pi}{2}\right).
\]

For $a_1 = \frac{1}{\sqrt{2}}, a_2 = \frac{1}{\sqrt{3}}$ :

\[\tag{4.6} \label{eqb13}
\sum\limits_{n=1}^{\infty} \left\{ \frac{\cot(\pi\sqrt{n^2+\frac{1}{6}})}{(n^2+\frac{1}{2})\sqrt{n^2+\frac{1}{6}}}+\frac{\cot(\pi\sqrt{n^2-\frac{1}{6}})}{(n^2+\frac{1}{3})\sqrt{n^2-\frac{1}{6}}}\right\} = \frac{3\sqrt{6}}{2}\coth\left(\frac{\sqrt{6}\pi}{6}\right) - \frac{\sqrt{6}\pi}{2}\coth\left(\frac{\sqrt{2}\pi}{2}\right)\coth\left(\frac{\sqrt{3}\pi}{3}\right)-\sqrt{6}\cot\left(\frac{\sqrt{6}\pi}{6}\right) .
\]

For $t=0, x_1 = 1, x_2 =i$, \eqref{eqt5} gives : 
\vspace{0.7cm}
\begin{corollary}
        For complex values $a_1,a_2$ such that $a_2^2-a_1^2$ is not a sum of two squares, there holds :
        \[ \tag{4.7} \label{eqt7}
        \sum\limits_{n=1}^{\infty} \frac{(2n^2+a_1^2-a_2^2)\coth(\pi \sqrt{n^2+a_1^2-a_2^2})}{(n^2+a_1^2)(n^2-a_2^2)\sqrt{n^2+a_1^2-a_2^2}} = \frac{\sqrt{a_1^2-a_2^2}}{2a_1^2a_2^2}\coth(\pi \sqrt{a_1^2-a_2^2}) - \frac{\pi}{2a_1a_2}\coth(\pi a_1)\cot(\pi a_2).
        \]
\end{corollary}

We shall point out the following case. For $a_1 = \frac{1}{\sqrt{2}}, a_2 = \frac{1}{\sqrt{3}}$ :

\[\tag{4.8} \label{eqb14}
\sum\limits_{n=1}^{\infty} \frac{(2n^2+\frac{1}{6})\coth(\pi \sqrt{n^2+\frac{1}{6}})}{(n^2+\frac{1}{2})(n^2-\frac{1}{3})\sqrt{n^2+\frac{1}{6}}} = \frac{\sqrt{6}}{2}\coth(\frac{\pi \sqrt{6}}{6}) - \frac{\pi \sqrt{6}}{2} \coth(\frac{\pi \sqrt{2}}{2})\cot(\frac{\pi \sqrt{3}}{3}).
\]

Now again in \eqref{eqt5} put $a_1 \to 0, a_2 \to 0$ and $t=1$ to arrive at:
\[
    \tag{4.9} \label{eqb15}
1-\frac{\pi^2}{3}(\frac{1}{x_1^2}+
\frac{1}{x_2^2})-\frac{2\pi}{x_1x_2}\sum\limits_{n=1}^{\infty}\left\{ \frac{\cot(\frac{n \pi x_1}{x_2})}{n(x_1^2n^2-1)} + \frac{\cot(\frac{n \pi x_2}{x_1})}{n(x_2^2n^2-1)} \right\} = \frac{\pi^2}{x_1x_2}\cot{\left(\ \frac{\pi}{x_1}\right)}\cot{\left(\ \frac{\pi}{x_2}\right)}.
\]

Making the transformations $x_1 \to \frac{1}{ix_1}$, $x_2 \to \frac{1}{x_2}$, where $i=\sqrt{-1}$, we arrive at :

\vspace{0.7cm}
\begin{corollary} For $x_1$ and $x_2$ such that $x_1^2n^2+x_2^2k^2 \neq 0$ for all natural numbers $k$ and $n$ :

\[
    \tag{4.10} \label{eqt8}
\pi^2x_1x_2\coth{(\pi x_1)}\cot{(\pi x_2)}= 1+ \frac{\pi^2}{3}(x_1^2-x_2^2) - 2\pi x_1x_2\left\{ x_1^2 \sum\limits_{n=1}^{\infty}\frac{\coth{\left(\frac{n\pi x_2}{x_1}\right)}}{n(n^2+x_1^2)}+x_2^2 \sum\limits_{n=1}^{\infty}\frac{\coth{\left(\frac{n\pi x_1}{x_2}\right)}}{n(n^2-x_2^2)}\right\}.
\]

\end{corollary}
R.Sitaramachandrarao in his work [2] established the above identity. Then using the elementary identities : 
\[
\frac{x^2}{n(n^2+x^2)}=-\frac{n}{n^2+x^2}+\frac{1}{n},
\]
\[
\frac{x^2}{n(n^2-x^2)}=\frac{n}{n^2-x^2}-\frac{1}{n},
\]
R.Sitaramachandrarao was able to arrive at : 
 \begin{multline}
    \tag{4.11} \label{eqb16}
\pi^2x_1x_2\coth{(\pi x_1)}\cot{(\pi x_2)}= 1+ \frac{\pi^2}{3}(x_1^2-x_2^2) + 2\pi x_1x_2\left\{  \sum\limits_{n=1}^{\infty}\frac{n\coth{\left(\frac{n\pi x_2}{x_1}\right)}}{n^2+x_1^2}- \sum\limits_{n=1}^{\infty}\frac{n\coth{\left(\frac{n\pi x_1}{x_2}\right)}}{n^2-x_2^2}\right\} \\- 4\pi x_1 x_2 \sum\limits_{n=1}^{\infty} \frac{1}{n}\left( \frac{1}{e^{2\pi n \frac{x_2}{x_1}} - 1} - \frac{1}{e^{2\pi n \frac{x_1}{x_2}} - 1}\right).
\end{multline}
Using this and the transformation formula for the logarithm of the Dedekind eta function, Bruce C. Berndt and Armin Straub [3] provided the following corrected version of \textbf{Entry 19(i)} by Ramanujan [1, p. 271]. For arbitrary $w$ and positive numbers $\alpha$ and $\beta$ such that $\alpha\beta = \pi^2$ :

\[ \tag{4.12} \label{eqb17}
\frac{\pi}{2}\cot(\sqrt{w\alpha})\coth(\sqrt{w\beta}) = \frac{1}{2w}+\frac{1}{2}\log\frac{\beta}{\alpha}+\sum\limits_{n=1}^{\infty}\left\{ \frac{n \alpha \coth(n \alpha)}{w + n^2\alpha}+\frac{n \beta \coth(n \beta)}{w - n^2\beta}\right\}.
\]

Let us recall the following Taylor series expansions :

\[ \tag{4.13} \label{eqb18}
x\coth(x) = \sum\limits_{k=0}^{\infty}\frac{B_{2k}(2x)^{2k}}{(2k)!}, 
\]

\[ \tag{4.14} \label{eqb19}
x\cot(x) = \sum\limits_{k=0}^{\infty}\frac{(-1)^kB_{2k}(2x)^{2k}}{(2k)!},
\]

where $B_{2k}$ are the Bernoulli numbers. Now in \eqref{eqt8}, let $x_1 \to x_1t$ and $x_2 \to x_2t$. Expanding both sides as power series of $t$ with the above Taylor expansions in mind and comparing the coefficients of $t^{2m+2}, m \geq 1$ we arrive at:

\[ \tag{4.15} \label{eqb20}
    2^{2m+1}\pi^{2m+1}\sum\limits_{k=0}^{m+1} \frac{{(-1)^k} B_{2k} B_{2m+2-2k} {x_2}^{2k} {x_1}^{2m+2-2k}}{(2k)!(2m+2-2k)!} =(-1)^{m}{x_2} {x_1}^{2m+1}\sum\limits_{n=1}^{\infty} \frac{\coth(\frac{n\pi {x_2}}{x_1})}{n^{2m+1}}-{x_1} {x_2}^{2m+1}\sum\limits_{n=1}^{\infty} \frac{\coth(\frac{n\pi {x_1}}{x_2})}{n^{2m+1}}.
\]  
Which is equivalent to Ramanujan's formula for $\zeta(2n+1)$ given in [1, p. 275-276].
\vspace{0.7cm}
\\Next, let us investigate further on \textbf{Proposition 2}. Put $M=2$ in \textbf{Proposition 2} to get the following statement.
\vspace{0.7cm}
\begin{corollary}
For $|\theta_1|, |\theta_2| \leq \frac{\pi}{2}$, and arbitrary $t$, and for complex values $a_1,a_2$ and nonzero complex values $x_1,x_2$ such that for all odd numbers $n$ and $k$, we have $x_1^2n^2-x_2^2k^2 \neq a_2^2-a_1^2$ and if $\textbf{Im}(\frac{x_1}{x_2}) \neq 0$ or $x_1 = x_2$ , there holds :
    \begin{multline} \tag{4.16} \label{eqt9} 
 x_1\sum\limits_{n=0}^{\infty}
\frac{(-1)^n\sin((2n+1)\theta_1)\sin\left(\frac{\theta_2 \sqrt{x_1^2(2n+1)^2+a_1^2-a_2^2}}{x_2}\right)}{(x_1^2(2n+1)^2+a_1^2-t^2)\sqrt{x_1^2(2n+1)^2+a_1^2-a_2^2}\cos\left(\frac{\pi \sqrt{x_1^2(2n+1)^2+a_1^2-a_2^2}}{2x_2}\right)}\\+ x_2\sum\limits_{n=0}^{\infty}
\frac{(-1)^n\sin((2n+1)\theta_2)\sin\left(\frac{\theta_1 \sqrt{x_2^2(2n+1)^2+a_2^2-a_1^2}}{x_1}\right)}{(x_2^2(2n+1)^2+a_2^2-t^2)\sqrt{x_2^2(2n+1)^2+a_2^2-a_1^2}\cos\left(\frac{\pi \sqrt{x_2^2(2n+1)^2+a_2^2-a_1^2}}{2x_1}\right)}\\ =\frac{\pi \sin\left(\frac{\theta_1 \sqrt{t^2-a_1^2}}{x_1}\right)\sin\left(\frac{\theta_2 \sqrt{t^2-a_2^2}}{x_2}\right)}{4\sqrt{t^2-a_1^2}\sqrt{t^2-a_2^2}\cos\left(\frac{\pi \sqrt{t^2-a_1^2}}{2x_1}\right)\cos\left(\frac{\pi \sqrt{t^2-a_2^2}}{2x_2}\right)}.
\end{multline}
\end{corollary}
Let $a_1=a_2=0$, then we arrive at :
\begin{multline} \tag{4.17} \label{eqb21} 
 \sum\limits_{n=0}^{\infty}
\frac{(-1)^n\sin((2n+1)\theta_1)\sin\left(\frac{(2n+1)\theta_2x_1}{x_2}\right)}{(2n+1)(x_1^2(2n+1)^2-t^2)\cos\left(\frac{(2n+1)\pi x_1}{2 x_2}\right)}+ \sum\limits_{n=0}^{\infty}
\frac{(-1)^n\sin((2n+1)\theta_2)\sin\left(\frac{(2n+1)\theta_1x_2}{x_1}\right)}{(2n+1)(x_2^2(2n+1)^2-t^2)\cos\left(\frac{(2n+1)\pi x_2}{2 x_1}\right)}\\ =\frac{\pi \sin\left(\frac{\theta_1 t}{x_1}\right)\sin\left(\frac{\theta_2 t}{x_2}\right)}{4t^2\cos\left(\frac{\pi t}{2x_1}\right)\cos\left(\frac{\pi t}{2x_2}\right)}.
\end{multline}

Make the transformations $x_1 \to \frac{1}{ x_1}$ and $x_2 \to \frac{1}{i x_2}$, where $i = \sqrt{-1}$ then we arrive at the following statement.

\vspace{0.7cm}

\begin{corollary}  For $|\theta_1|, |\theta_2| \leq \frac{\pi}{2}$, and arbitrary $t$ and for $x_1$ and $x_2$ such that $x_1^2n^2+x_2^2k^2 \neq 0 $ for all odd numbers $n$ and $k$, there holds :

\begin{multline}
    \frac{\pi \sin\left(\theta_1x_1 t\right)\sinh\left(\theta_2x_2 t\right)}{4t^2\cos\left(\frac{\pi x_1 t}{2}\right)\cosh\left(\frac{\pi x_2 t}{2}\right)} = x_2^2\sum\limits_{n=0}^{\infty}
\frac{(-1)^n\sin((2n+1)\theta_2)\sinh\left(\frac{(2n+1)\theta_1x_1}{x_2}\right)}{(2n+1)((2n+1)^2+t^2x_2^2)\cosh\left(\frac{(2n+1)\pi x_1}{2 x_2}\right)}+ x_1^2\sum\limits_{n=0}^{\infty}
\frac{(-1)^n\sin((2n+1)\theta_1)\sinh\left(\frac{(2n+1)\theta_2x_2}{x_1}\right)}{(2n+1)((2n+1)^2-t^2x_1^2)\cosh\left(\frac{(2n+1)\pi x_2}{2x_1}\right)}.\\ \tag{4.18} \label{eqt10}
\end{multline}
\end{corollary}We see that this is in fact identical to \textbf{Corollary 2} by Ramanujan [1, p. 269]. Now we proceed to examine \eqref{eqt9}. Divide both sides of the equality by $\theta_1 \theta_2$ and let $\theta_1 \to 0$, $\theta_2 \to 0$, put $x_1=1,x_2=i, t=0$, then we arrive to the following statement.
\vspace{0.7cm}
\begin{corollary} For $a_1$ and $a_2$ such that $a_2^2-a_1^2$ is not a sum of two odd squares, there holds :
    \[\tag{4.19} \label{eqt11}
    \sum\limits_{n=0}^{\infty} \frac{(-1)^n(2n+1)(2(2n+1)^2+a_1^2-a_2^2)}{((2n+1)^2+a_1^2)((2n+1)^2-a_2^2)\cosh\left( \frac{\pi}{2} \sqrt{(2n+1)^2+a_1^2-a_2^2}\right)} = \frac{\pi}{4\cosh\left(\frac{\pi a_1}{2} \right)\cos\left( \frac{\pi a_2}{2}\right)}.
    \]
\end{corollary}

Put $a_1 = \frac{1}{\sqrt{2}}, a_2= \frac{1}{\sqrt{3}}$ to get :
\[\tag{4.20} \label{eqb22}
    \sum\limits_{n=0}^{\infty} \frac{(-1)^n(2n+1)(2(2n+1)^2+\frac{1}{6})e^{\frac{\pi}{2} \sqrt{(2n+1)^2+\frac{1}{6}}}}{((2n+1)^2+\frac{1}{2})((2n+1)^2-\frac{1}{3})\left( e^{\pi \sqrt{(2n+1)^2+\frac{1}{6}}}+1\right)} = \frac{\pi}{8\cosh\left(\frac{\pi \sqrt{2}}{4} \right)\cos\left( \frac{\pi \sqrt{3}}{6}\right)}.
    \]

In \textbf{Proposition 2}, divide both sides of the equality by $\theta_1 \theta_2$ and let $\theta_1 \to 0$, $\theta_2 \to 0$, then we arrive at the following statement. 

\vspace{0.7cm}

\begin{corollary} 
For arbitrary $t$, and for complex values $a_1,...,a_M$ and nonzero complex values $x_1,...,x_M$ such that for all $i \neq j$, and for all odd natural numbers $n$ and $k$, we have $x_i^2n^2-x_j^2k^2 \neq a_j^2-a_i^2$, and if for each $i \neq j$, $\textbf{Im}(\frac{x_i}{x_j}) \neq 0$, there holds :
    \[\tag{4.21} \label{eqt12}
\sum\limits_{1 \leq i \leq M} x_i^2\sum\limits_{n=0}^{\infty}
\frac{(-1)^n(2n+1)}{x_i^2(2n+1)^2+a_i^2-t^2}\prod\limits_{\substack{1 \leq j \leq M \\ j \neq i }}\frac{1}{\cos\left(\frac{\pi \sqrt{x_i^2(2n+1)^2+a_i^2-a_j^2}}{2x_j}\right)} = \frac{\pi}{4}\prod\limits_{1\leq i \leq M}\frac{1}{\cos\left(\frac{\pi \sqrt{t^2-a_i^2}}{2x_i}\right)}.
\]
\end{corollary}

We point out the particular case for $M=3, a_1=a_2=a_3=0$ in \eqref{eqt12}, to arrive at the following result. 
\vspace{0.7cm}
\begin{corollary} For arbitrary $t$ and for complex values $x_1,x_2,x_3$ such that $\textbf{Im}(\frac{x_1}{x_2}) \neq 0$, $\textbf{Im}(\frac{x_3}{x_2}) \neq 0$, and $\textbf{Im}(\frac{x_1}{x_3}) \neq 0$, there holds :
   \begin{multline}   \tag{4.22} \label{eqt13}
    x_1^2\sum\limits_{n=0}^{\infty}
\frac{(-1)^n(2n+1)}{(x_1^2(2n+1)^2-t^2)\cos\left(\frac{(2n+1)\pi x_1}{2x_2}\right)\cos\left(\frac{(2n+1)\pi x_1}{2x_3}\right)}+x_2^2\sum\limits_{n=0}^{\infty}
\frac{(-1)^n(2n+1)}{(x_2^2(2n+1)^2-t^2)\cos\left(\frac{(2n+1)\pi x_2}{2x_1}\right)\cos\left(\frac{(2n+1)\pi x_2}{2x_3}\right)}\\+x_3^2\sum\limits_{n=0}^{\infty}
\frac{(-1)^n(2n+1)}{(x_3^2(2n+1)^2-t^2)\cos\left(\frac{(2n+1)\pi x_3}{2x_1}\right)\cos\left(\frac{(2n+1)\pi x_3}{2x_2}\right)} = \frac{\pi}{4 \cos\left(\frac{\pi t}{2x_1}\right)\cos\left(\frac{\pi t}{2x_2}\right)\cos\left(\frac{\pi t}{2x_3}\right)}.
\end{multline} 
\end{corollary}

We now state the following formula for $\frac{\pi}{4}$. Let $a_1 = ... = a_M = 0$ and let $t = 0$ in \eqref{eqt12} to arrive at the following result.
\vspace{0.7cm}
\begin{corollary} For complex values $x_1,...,x_M$ such that for all $i \neq j$, $\textbf{Im}(\frac{x_i}{x_j}) \neq 0$, there holds :
    \[\tag{4.23} \label{eqt14}
\frac{\pi}{4} = \sum\limits_{1 \leq i \leq M} \sum\limits_{n=0}^{\infty}
\frac{(-1)^n}{2n+1}\prod\limits_{\substack{1 \leq j \leq M \\ j \neq i }}\frac{1}{\cos\left(\frac{(2n+1)\pi x_i}{2x_j}\right)}.
\]
\end{corollary}

\section{On a general series}

In [6] and [7], Dixit and Maji, and Kanemitsu, Tanigawa and Yoshimoto respectively studied the generalized Lambert series of the form : \textit{\begin{equation}\notag
\frac{1^{r}}{e^{1^sx}-1}+\frac{2^{r}}{e^{2^sx}-1}+\frac{3^{r}}{e^{3^sx}-1}+\cdots
\end{equation}
where $s$ is a positive integer and $r-s$ is any even integer.} 
\\It is also noted that Ramanujan mentioned this series in page 332 of the Lost Notebook [8]. Motivated by these papers, using our method of partial fractions already described by the general relations \eqref{eqt3} and \eqref{eqab10} we shall now here try to reach towards a formula regarding to the more general series : 

\[
\sum\limits_{n=1}^{\infty}\frac{1}{n^{2s-\frac{s}{r}}(n^{2s}-z)(e^{n^{\frac{s}{r}}x}-1)}
\]

\emph{where $s$ and $r$ are natural numbers.}
\\For our purpose we shall use the following reduced form of \eqref{eqt3}. Let $f_1$ and $f_2$ be two functions having expressions of the form :
\[
f_i(z) = \sum\limits_{n=1}^{\infty}\frac{R_{i}(n)}{c_{i}(n)-z}, i= 1,2.
\]
Then, for arbitrary $t$ and for complex values $x_1,x_2$ such that $x_1c_1(n)-x_2c_2(k) \neq 0$ for all natural numbers $n$ and $k$, there holds :
    \[\tag{5.1} \label{eqs4basicformula}
        f_1\left(\frac{t}{x_1}\right) f_2\left(\frac{t}{x_2}\right) =  x_1\sum\limits_{n=1}^{\infty}\frac{R_1(n)}{x_1c_1(n)-t} f_2\left(\frac{x_1c_1(n)}{x_2}\right)+x_2\sum\limits_{n=1}^{\infty}\frac{R_2(n)}{x_2c_2(n)-t} f_1\left(\frac{x_2c_2(n)}{x_1}\right).
    \]
Now let :

\[f_1(z) = \sum\limits_{n=1}^{\infty} \frac{1}{n^{2M}-z}, \hspace{0.05cm} f_2(z) = \sum\limits_{n=1}^{\infty} \frac{1}{n^{2N}-z}, \text{where $M$ and $N$ are natural numbers}.\]

Since it holds that :

\[ \tag{5.2} \frac{1}{n^{2N} - z} = \frac{ z^{\frac{1}{N}-1}}{N} \sum\limits_{0 \leq r \leq N-1} \frac{e^{\frac{2ir\pi}{N}}}{n^2-z^{\frac{1}{N}}e^{\frac{2ir\pi}{N}}}. \]

Thus we have :

\[f_1(z) = \frac{1}{2z}-\frac{\pi}{2M}\sum\limits_{0 \leq r \leq M-1} \frac{e^{\frac{ir\pi}{M}}}{z^{1-\frac{1}{2M}}}\cot(\pi z^{\frac{1}{2M}}e^{\frac{ir\pi}{M}}), \hspace{0.05cm} f_2(z) =  \frac{1}{2z} - \frac{\pi}{2N}\sum\limits_{0 \leq r \leq N-1} \frac{e^{\frac{ir\pi}{N}}}{z^{1-\frac{1}{2N}}} \cot(\pi z^{\frac{1}{2N}}e^{\frac{ir\pi}{N}}).\]

Therefore via \eqref{eqs4basicformula}, 

\begin{multline}
    \left(\frac{x_1}{2t}-\frac{\pi x_1^{1-\frac{1}{2M}}}{2Mt^{1-\frac{1}{2M}}}\sum\limits_{0 \leq r \leq M-1} e^{\frac{ir\pi}{M}} \cot\left(\pi \left(\frac{t}{x_1}\right)^{\frac{1}{2M}}e^{\frac{ir\pi}{M}}\right)\right)\left(\frac{x_2}{2t} - \frac{\pi x_2^{1-\frac{1}{2N}}}{2N t^{1-\frac{1}{2N}}}\sum\limits_{0 \leq r \leq N-1} e^{\frac{ir\pi}{N}} \cot\left(\pi \left( \frac{t}{x_2}\right)^{\frac{1}{2N}}e^{\frac{ir\pi}{N}}\right)\right) = \frac{x_2}{2x_1} \sum\limits_{n=1}^{\infty} \frac{1}{( n^{2M}-\frac{t}{x_1})n^{2M}} + \frac{x_1 }{2 x_2} \sum\limits_{n=1}^{\infty} \frac{1}{( n^{2N}-\frac{t}{x_2}) n^{2N}} \\ \tag{5.3} - \frac{\pi x_1^{\frac{1}{2N}}x_2^{1-\frac{1}{2N}}}{2N}\sum\limits_{0 \leq r \leq N-1} e^{\frac{ir\pi}{N}} \sum\limits_{n=1}^{\infty} \frac{\cot(\pi e^{\frac{ir\pi}{N}} (\frac{x_1}{x_2})^{\frac{1}{2N}} n^{\frac{M}{N}})}{n^{2M-\frac{M}{N}}(x_1n^{2M}-t)} - \frac{\pi x_2^{\frac{1}{2M}}x_1^{1-\frac{1}{2M}}}{2M}\sum\limits_{0 \leq r \leq M-1} e^{\frac{ir\pi}{M}} \sum\limits_{n=1}^{\infty} \frac{\cot(\pi e^{\frac{ir\pi}{M}} (\frac{x_2}{x_1})^{\frac{1}{2M}} n^{\frac{N}{M}})}{n^{2N-\frac{N}{M}}(x_2n^{2N}-t)}.
\end{multline}

Now we put $x_1 \to x_1^{2NM}$, $x_2 \to x_2^{2NM}$, $t \to t^{2NM}$ and in view of the Eulerian formula $ \zeta(2n) = \frac{(-1)^{n-1}(2\pi)^{2n}B_{2n}}{2(2n)!}$, where $B_{2n}$ are Bernoulli numbers [9], and making further manipulations and substitutions we arrive at :

\begin{multline}
    \frac{\pi x_1^{M-2MN}}{Nx_2^{M-2MN}} \sum\limits_{0 \leq r \leq N-1} e^{\frac{ir\pi}{N}} \sum\limits_{n=1}^{\infty} \frac{\cot\left( \pi e^{\frac{ir\pi}{N}}\frac{x_1^M}{x_2^M} n^{\frac{M}{N}}\right)}{n^{2M-\frac{M}{N}}(n^{2M} - \frac{t^{2MN}}{x_1^{2MN}})} + \frac{\pi x_2^{N-2MN}}{Mx_1^{N-2MN}} \sum\limits_{0 \leq r \leq M-1} e^{\frac{ir\pi}{M}} \sum\limits_{n=1}^{\infty} \frac{\cot\left( \pi e^{\frac{ir\pi}{M}}\frac{x_2^N}{x_1^N} n^{\frac{N}{M}}\right)}{n^{2N-\frac{N}{M}}(n^{2N} - \frac{t^{2MN}}{x_2^{2MN}})} = \frac{x_1^{2MN}x_2^{2MN}}{2t^{4MN}}+ \frac{x_2^{2MN} (-1)^{M}(2\pi)^{2M}B_{2M}}{2t^{2MN}(2M)!}\\ \tag{5.4}+\frac{x_1^{2MN} (-1)^{N}(2\pi)^{2N}B_{2N}}{2t^{2MN}(2N)!}-\frac{\pi^2 x_1^{2MN-N}x_2^{2MN-M}}{2MN t^{4MN-N-M}}\left\{ \sum\limits_{0 \leq r \leq M-1} e^{\frac{ir\pi}{M}}\cot \left( \frac{\pi t^{N} e^{\frac{ir\pi}{M}}}{x_1^N}\right) \right\}\left\{ \sum\limits_{0 \leq r \leq N-1} e^{\frac{ir\pi}{N}}\cot \left( \frac{\pi t^{M} e^{\frac{ir\pi}{N}}}{x_2^{M}}\right) \right\}.
\end{multline}

Now we put $M \geq N$ and let $x_1 \to e^{\frac{i \pi}{2N}}x_1$, using the fact $\coth(x) = 1+ \frac{2}{e^{2x}-1}$, we finally arrive at the following proposition. 
\vspace{0.7cm}
\begin{theorem} For natural numbers $M$ and $N$ with $M \geq N$, and arbitrary $t$, and for complex values $x_1$ and $x_2$ such that $x_1^{2MN}n^{2M}+(-1)^{M-1}x_2^{2MN}k^{2N} \neq 0$ for all natural numbers $n$ and $k$, there holds :
\begin{multline}
\frac{\pi x_2^{N-2MN}}{Mx_1^{N-2MN}(1-e^{\frac{i\pi}{M}})} \sum\limits_{n=1}^{\infty}\frac{1}{n^{2N-\frac{N}{M}}(n^{2N}-\frac{t^{2MN}}{x_2^{2MN}})}+\frac{\pi x_2^{N-2MN}}{Mx_1^{N-2MN}} \sum\limits_{0 \leq r \leq M-1}  e^{\frac{ir\pi}{M}} \sum\limits_{n=1}^{\infty} \frac{1}{n^{2N-\frac{N}{M}}(n^{2N}-\frac{t^{2MN}}{x_2^{2MN}})\left( e^{\frac{2n^{\frac{N}{M}}\pi e^{\frac{ir\pi}{M}}x_2^N}{x_1^N}} -1\right)} \\ +e^{\frac{i\pi(M-N)}{2N}}\frac{\pi x_1^{M-2MN}}{Nx_2^{M-2MN}(1-e^{\frac{i\pi}{N}})}\sum\limits_{n=1}^{\infty}\frac{1}{n^{2M-\frac{M}{N}}(n^{2M}+(-1)^{M-1}\frac{t^{2MN}}{x_1^{2MN}})}+ e^{\frac{i\pi(M-N)}{2N}}\frac{\pi x_1^{M-2MN}}{Nx_2^{M-2MN}}\sum\limits_{0 \leq r \leq N-1} e^{\frac{ir\pi}{N}} \sum\limits_{n=1}^{\infty}\frac{1}{n^{2M-\frac{M}{N}}(n^{2M}+(-1)^{M-1}\frac{t^{2MN}}{x_1^{2MN}})\left(e^{\frac{2n^{\frac{M}{N}}\pi e^{\frac{i\pi(r+\frac{M}{2}-\frac{N}{2})}{N}}x_1^M}{x_2^M}}-1\right)}\\ \tag{5.5} = \frac{x_1^{2MN}x_2^{2MN}}{4t^{4MN}} + \frac{x_2^{2MN}(2\pi)^{2M}B_{2M}}{4t^{2MN}(2M)!}+\frac{x_1^{2MN}(-1)^{N}(2\pi)^{2N}B_{2N}}{4t^{2MN}(2N)!}-\frac{\pi^2 x_1^{2MN-N}x_2^{2MN-M}}{4MNt^{4MN-N-M}}\left\{ \sum\limits_{0 \leq r \leq M-1} e^{\frac{ir\pi}{M}} \coth\left(\frac{\pi t^{N} e^{\frac{ir\pi}{M}}}{x_1^{N}} \right) \right\} \left\{ \sum\limits_{0 \leq r \leq N-1} e^{\frac{ir\pi}{N}} \cot\left( \frac{\pi t^M e^{\frac{ir\pi }{N} }}{x_2^M}\right) \right\}.
\end{multline}

Provided that the series from the left hand side converges absolutely.
\end{theorem}

\end{document}